\documentclass[12pt]{article}
\usepackage{amssymb,latexsym,amsfonts,amsbsy,color,setspace}

\newcommand{\qed}{\hbox{\rule[-2pt]
{3pt}{6pt}}}

\newtheorem{dfe}{Definition}
[section]
\newtheorem{lem}[dfe]{Lemma}
\newtheorem{theo}[dfe]{Theorem}

\newtheorem{pro}[dfe]{Proposition}
\newtheorem{cor}[dfe]{Corollary}

\makeatletter

\@addtoreset{equation}{section}
\makeatother


\usepackage{enumerate}
\usepackage{latexsym,amsmath,amsxtra,amsfonts,amsbsy,curves,bezier,epic,color,eepic,bezier,amscd,makeidx}
\newcommand{\ds}{\displaystyle}

\begin{title}
{Relative $t$-designs in binary Hamming association scheme $H(n,2)$}
\end{title}
\author{Eiichi Bannai, Etsuko Bannai and Yan Zhu\\
{}\\
\normalsize\sl{
Dedicated to Andries Brouwer on the occasion of his 65th birthday }
}

\date{}
\begin{document}
\maketitle
\begin{abstract}
A relative $t$-design in the binary Hamming association schemes $H(n,2)$ is equivalent to a weighted regular $t$-wise balanced design,
i.e., certain combinatorial $t$-design which allow different sizes of blocks and a weight function on blocks.
In this paper, we study relative $t$-designs in $H(n,2)$, putting emphasis on Fisher type inequalities and the existence of tight relative $t$-designs.
We mostly consider relative $t$-designs on two shells.
We prove that if the weight function is constant on each shell of a relative $t$-design on two shells then the subset in each shell must be a combinatorial $(t-1)$-design.
This is a generalization of the result of Kageyama who proved this under the stronger assumption that the weight function is constant on the whole block set.  Using this, we define tight relative $t$-designs for odd $t$, and a strong restriction on the possible parameters of tight relative $t$-designs in $H(n,2)$.
We obtained a new family of such tight relative $t$-designs, which were unnoticed before.
We will give a list of feasible parameters of such relative 3-designs with $n \leq 100$, and then we discuss the existence and/or the non-existence of such tight relative $3$-designs.
We also discuss feasible parameters of tight relative 4-designs on two shells in $H(n,2)$ with $n \leq 50$.
In this study we come up with the connection on the topics of classical design theory, such as symmetric 2-designs (in particular $2$-$(4u-1,2u-1,u-1)$ Hadamard designs) and Driessen's result on the non-existence of certain 3-designs.
We believe the Problem 1 and Problem 2 presented in Section \ref{subsec:further_result} open a new way to study relative $t$-designs in $H(n,2)$.
We conclude our paper listing several open problems.
\end{abstract}


\section{Introduction}
The concept of relative $t$-designs in association schemes was started by Delsarte in \cite{Delsarte-1977} (1977).
 We refer to \cite{Bannai-B-Z-2015} for a survey of a history of the study of relative $t$-designs in association schemes.
In this paper we study relative $t$-designs in binary Hamming association scheme $H(n,2)$.
Our emphasis is on Fisher type inequalities and tight relative $t$-designs.
We mostly consider relative $t$-designs on two shells.

Two types of relative $t$-designs in any P-and Q-polynomial association schemes were considered, see \cite{B-B-S-Tanaka}.
In the case of Hamming association schemes $H(n,q)$, they coincide.
Firstly, we study the equivalence of these definitions, as well as the equivalence with the concept of weighted regular $t$-wise balanced designs (cf. \cite{Kramer-K-2006}).
The second topic is on the theory of relative $t$-designs on two shells.
Our main theorem is Theorem \ref{theo:3.3}, which says that if the weight function is constant on each shell of a relative $t$-design on two shells then the subset in each shell must be a combinatorial $(t-1)$-design.
This is a generalization of Proposition 1 of Kageyama in \cite{Kageyama-1991}.
Kageyama proved his Proposition 1 for $t$-wise and $(t-1)$-wise balanced design on two shells in $H(n,2)$, under the stronger assumption that the weight function is constant on the whole set.
Our theorem (Theorem \ref{theo:3.3}) gives a strong restriction on the possible parameters of tight relative $t$-designs in $H(n,2)$.
We studied tight relative 2-designs on two shells in Johnson association schemes $J(v,k)$ in \cite{Zhu-B-B-2015} (see also \cite{Bannai-B-Z-2015}).
The concept of tight relative $2e$-designs were already discussed in \cite{Li-B-B-2013}, \cite{Bannai-B-B-2014}, etc.
Here we first discuss Fisher type lower bound, and define tight relative 3-designs on two shells in $H(n,2)$.
 We discuss some examples, as well as their classification problems.
 In Section \ref{sec:design_odd}, we will give a list of feasible parameters of such relative 3-designs with $n\leq 100$.
 Then we discuss the existence and/or the non-existence of such tight relative $3$-designs.
In Section \ref{sec:4-design}, we discuss possible feasible parameters of tight relative 4-designs in $H(n,2)$ with $n\leq 50$.
(Here note that our theorems in the previous sections play a very important role.
We discuss the existence (and mostly non-existence) results of tight relative 4-designs on two shells in $H(n,2)$ with $n\leq 50$.
 We conclude this paper by mentioning further research problems, namely what are the problems we want to study in this research direction.

\section{Definitions and basic facts}
Let $\mathfrak X=(X,\{R_i\}_{0\leq i\leq n})$ be a symmetric association scheme defined on $X$. (Please refer to \cite{Delsarte-1973}, \cite{Bannai-I-1984}, \cite{Brouwer-C-N-1989} for information on association schemes.)
Let $\mathcal F(X)$ be the vector space of all the real valued functions defined on $X$.
Let $u_0\in X$ be fixed.
Let $X_j=\{x\in X\mid (x,u_0)\in R_j\}$.
Some designs in symmetric association schemes are defined very similar to spherical designs and Euclidean designs.
When we consider spherical designs or Euclidean designs, we use the vector space of polynomials.
When we try to define designs in symmetric association schemes by similar manner as spherical or Euclidean designs, we use the space $\mathcal F(X)$ of real valued functions on $X$ instead of the space of polynomials.
We consider $\mathcal F(X)$ in two different ways.
One way is to study $\mathcal F(X)$ using the property of P-structure of $\mathfrak X$.
For $z\in X_j$ we define a real valued function $f_z$ on $X$ in the following way.
$$f_z(x)=\left\{
\begin{array}{ll}
1&\mbox{if $x\in X_i$, $i\geq j$ and $(x,z)\in R_{i-j}$},\\
0&\mbox{otherwise.}
\end{array}
\right.$$
We define
$$\mathrm{Hom}_j(X)=\langle f_z\mid z\in X_j\rangle.$$
Let $k_j$ be the $j$-th valency of $\mathfrak X$. Then $\dim(\mathrm{Hom}_j(X))=k_j$, and we have the following direct sum.
$$\mathcal F(X)=\mathrm{Hom}_0(X)+\mathrm{Hom}_1(X)+\cdots+\mathrm{Hom}_n(X).$$
Another way to study $\mathcal F(X)$ is using the column space of the primitive idempotents $E_0,E_1,\ldots,E_n$ of $\mathfrak X$.
For each $\ell,~ 0\leq \ell\leq n$, and $u\in X$, let $\phi^{(\ell)}_u$ be the $u$-th column vector of $|X|E_\ell$.
Each $\phi^{(\ell)}_u$ is also regarded as a function defined on $X$.
For each $E_\ell$ we define
$$
L_\ell(X)=\langle \phi^{(\ell)}_u \mid u\in X \rangle, \quad 0\leq \ell \leq n.
$$
Let $m_\ell$ be the rank of $E_\ell$.
Then $\dim(L_\ell(X))=m_\ell$.
If we consider the usual inner product of $\mathbb R^{|X|}$, then we have an orthogonal decomposition of $\mathcal F(X)$.
$$\mathcal F(X)=L_0(X)\perp L_1(X)\perp \cdots\perp L_n(X).$$
Let $\mathfrak X=(X,\{R_i\}_{0\leq i\leq n})$ be an association scheme defined on $X$.
Let $(Y,w)$ be a positive weighted subset of $X$.
Assume $Y$ is on a union of $p$ shells $X_{r_1}\cup\cdots\cup X_{r_p}$ in $\mathfrak X$.
Let $Y_{r_\nu}=Y\cap X_{r_\nu}$ and $w(Y_{r_\nu})=\sum_{y\in Y_{r_\nu}}w(y)$ for $1\leq \nu \leq p$.
We have the following two different ways of the definition for relative $t$-designs (cf. \cite{B-B-S-Tanaka}).
%%%%%%%
\begin{dfe}[Relative $t$-design in Q-polynomial association scheme]
$ $\\
Let $\mathfrak X$ be a Q-polynomial association scheme of class $n$.
Then $(Y,w)$ is a relative $t$-design of the Q-structure of $\mathfrak X$ with respect to $u_0$, if the following property holds:
\begin{eqnarray}
&&\sum_{\nu=1}^p\frac{w(Y_{r_\nu})}{|X_{r_\nu}|}\sum_{x\in X_{r_\nu}}f(x)
=\sum_{y\in Y}w(y)f(y)
\label{equ:2.1}
\end{eqnarray}
for any $f\in L_j(X)$ with $1\leq j\leq t$.
\end{dfe}
\begin{dfe}[Relative $t$-design in P-polynomial association scheme]
$ $\\
Let $\mathfrak X$ be a P-polynomial association scheme of class $n$.
Then $(Y,w)$ is a relative $t$-design of the P-structure of $\mathfrak X$ with respect to $u_0$, if the following property holds:
$$\sum_{\nu=1}^p\frac{w(Y_{r_\nu})}{|X_{r_\nu}|}\sum_{x\in X_{r_\nu}}f(x)=\sum_{y\in Y}w(y)f(y)$$
for any $f\in \mathrm{Hom}_j(X)$ with $1\leq j\leq t$.
\end{dfe}
As for the relative $2e$-designs in Q-polynomial association schemes, the following Fisher type lower bound is known.
\begin{theo}[\cite{Bannai-B-2012}]
Let $(Y,w)$ be a relative $2e$-design of a Q-polynomial association scheme with respect to $u_0$.
Assume $Y$ is on $p$ shells $S=X_{r_1}\cup X_{r_2}\cup\cdots\cup X_{r_p}$.
Then
\begin{eqnarray}
|Y|\geq \dim(L_e(S)+L_{e-1}(S)+\cdots+L_{e-p+1}(S))
\label{equ:2.2}
\end{eqnarray}
holds, where $L_j(S)$ denotes the restriction of $L_j(X)$ to $S$.
\label{theo:2.2}
\end{theo}
\begin{dfe}[\cite{Bannai-B-2012}]
If equality holds in \eqref{equ:2.2}, $(Y,w)$ is called tight.
\end{dfe}
Theorem \ref{theo:2.2} is proved by using the following proposition.
\begin{pro}[\cite{Bannai-B-2012}]
For $u_1,u_2\in X$, let $\phi^{(\ell_1)}_{u_1}$ and $\phi^{(\ell_2)}_{u_2}$ be the $u_1$-th column vector of $E_{\ell_1}$ and $u_2$-th column vector of $E_{\ell_2}$ respectively.
Then $(E_{\ell_1}\phi_{u_1})(E_{\ell_2}\phi_{u_2})\in \sum_{j=0}^{2\ell}L_j(X)$ holds for any integers $\ell_1, \ell_2, \ell\in \{0,1,\ldots,n\}$ satisfying $\ell_1+\ell_2\leq 2\ell$ .
\label{pro:2.5}
\end{pro}
Thus for a P-and Q-polynomial association scheme we have two types of decomposition of $\mathcal F(X)$.
In general for a P-polynomial association scheme $\mathrm{Hom}_{0}(X)+\cdots+\mathrm{Hom}_{\ell}(X)$ may not be closed under the product of functions.
However in \cite{B-B-S-Tanaka}, they proved that the following condition is satisfied for general Hamming association scheme $H(n,q)$.
$$
L_0(X)+\cdots+L_\ell(X)=\mathrm{Hom}_0(X)+\cdots+\mathrm{Hom}_\ell(X)\quad \mbox{for $\ell=0,\ldots,n$}.
$$
Proposition \ref{pro:2.5} implies the following lemma.
\begin{lem}
Let $(Y,w)$ be a positive weighted subset on $p$ shells $X_{r_1}\cup\cdots\cup X_{r_p}$ in $H(n,2)$.
Then $(Y,w)$ is a relative $t$-design with respect to $u_0$ if and only if the following condition holds:
\begin{eqnarray}
\sum_{\nu=1}^p\frac{w(Y_{r_\nu})}{|X_{r_i}|}\sum_{x\in X_{r_\nu}}\prod_{j=1}^s\phi^{(1)}_{u_j}(x)=\sum_{y\in Y}w(y)\prod_{j=1}^s\phi^{(1)}_{u_j}(y)
\label{equ:2.3}
\end{eqnarray}
for any $s$ with $1\leq s\leq t$ and any $u_1,\ldots, u_s\in X_1$.
\label{lem:2.6}
\end{lem}
{\bf Proof~} In \cite{Li-B-B-2013}, it is proved that $\{\phi^{(1)}_{u}\mid u\in X_1\}$ spans the column space of $E_1$.
Since $H(n,2)$ is a Q-polynomial association scheme, there exists a polynomial $v^*_i$ of degree $i$ and $|X|E_i=v^*_i(|X|E_1)$ holds.
Here product of the matrix $(|X|E_1)^j$ is defined by the Hadamard product.
Therefore if $Y\subset X_{r_1}\cup X_{r_2}\cup \cdots \cup X_{r_p}$ satisfies the condition \eqref{equ:2.3}, then $(Y,w)$ satisfies the defining equation \eqref{equ:2.1} of relative $t$-designs.
\hfill\qed\\

The equality \eqref{equ:2.3} plays an important role when we determine the feasible parameters of tight relative $t$-design in $H(n,2)$.

\section{Main results}\label{sec:main}

Let $V=\{1,2,\ldots,n\}$ and ${V\choose r}$ be the set of all $r$-point subsets of $V$.
Let $X=F_2^n$.
Without loss of generality we may assume $u_0=(0,\ldots, 0)\in X$ and $X_r=\{(x_1,\ldots,x_n)\in X\mid \sharp \{i\mid x_i=1\}=r\}$.
For $x=(x_1,\ldots,x_n)\in X$, we define a subset $B_x$ of $V$ by $B_x=\{i\in V\mid  x_i=1\}$.
For a subset $Y\subset X$, we define $\mathcal B_Y=\{B_y\mid y\in Y\}$.
On the other hand for a subset $B\in {V\choose r}$, we define $x_B=(x_1,\ldots,x_n)\in X_r$ by $x_i=1$ if $i\in B$ and $x_i=0$ if $i\not\in B$.
\begin{dfe}
Let $V=\{1,2,\ldots,n\}$ and $\mathcal B$ be a set of subsets in $V$. Let $w$ be a positive weight function on $\mathcal B$.
 Then $(V,\mathcal B,w)$ is a $j$-wise balanced design if
 $$\sum_{B\in \mathcal B, B_z\subset B}w(B)=\lambda_j$$
holds with a constant $\lambda_j$ which determined by $j$ and independent on the choice of $z\in X_j$.
$(V,\mathcal B,w)$ is called a regular $t$-wise balanced design, if it is $j$-wise balanced for $j=1,2,\ldots,t$.
\label{dfe:3.1}
\end{dfe}
\begin{theo}[\cite{Delsarte-S-1989}]
Let $V=\{1,2,\ldots,n\}$.
Let $Y\subset X$ and $\mathcal B_Y=\{B_y\mid y\in Y\}$.
Then $(Y,w)$ is a relative $t$-design of the P-structure of $H(n,2)$ on $p$ shells with respect to $u_0$, if and only if $(V,\mathcal B_Y,w)$ is a regular $t$-wise balanced design with positive weight $w$.
\label{theo:3.2}
\end{theo}
In the following for a subset $Y\subset X$, if $(V,\mathcal B_Y,w)$ is a regular $t$-wise balanced design with positive weight $w$, then we say $Y$ has the structure of a regular $t$-wise balanced design.
If $(V,\mathcal B_Y)$ is a combinatorial design, then we say $Y$ has the structure of a combinatorial design.
First we prove the following theorem, which is a generalization of Kageyama's Theorem \cite{Kageyama-1991}.
\begin{theo}[Generalization of Kageyama's Theorem]
$ $\\
Let $(V,\mathcal B,w)$ be a $t$-wise and $(t-1)$-wise balanced design.
Assume $\mathcal B$ consists of blocks of size $r_1$ and $r_2$.
Let $ \mathcal B_{r_1}=\{B\in \mathcal B\mid |B|=r_1\}$ and $ \mathcal B_{r_2}=\{B\in \mathcal B\mid |B|=r_2\}$.
Assume that the weight function $w$ is a constant $w_{r_\nu}$ on each block set $\mathcal B_{r_\nu}$, $\nu =1,2$.
Then $(V, \mathcal B_{r_\nu})$ is a combinatorial $(t-1)$-$(n,r_\nu, \lambda^{(r_\nu)}_{t-1})$ design for $\nu=1,2$, with
$$\begin{array}{l} \medskip
\ds \lambda^{(r_1)}_{t-1}=\frac{(r_2-t+1)\lambda_{t-1}-(n-t+1)\lambda_t}{(r_2-r_1)w_{r_1}},\\
\ds \lambda^{(r_2)}_{t-1}=\frac{(r_1-t+1)\lambda_{t-1}-(n-t+1)\lambda_t}{(r_1-r_2)w_{r_2}}.
\end{array}$$
Here $\displaystyle \lambda_{j}=\sum_{B\in \mathcal B, B_z\subset B}w(B)$ for $z\in X_{j}$, for $j=t-1,t$.
\label{theo:3.3}
\end{theo}
Then we prove the following theorem.
\begin{theo}
Definition and notation are given as above.
Let $(Y,w)$ be a weighted subset of $X_{r_1}\cup X_{r_2}$.
Let $Y_{r_\nu}=Y\cap X_{r_\nu}$ and assume $w(y)=w_{r_\nu}$ for any $y\in Y_{r_\nu}$ for $\nu=1,2$.
Let $\mathcal B_{r_\nu}$ be the set of blocks corresponding to $Y_{r_\nu}$  and $N_{r_\nu}=|Y_{r_\nu}|$ for $\nu=1,2$.
Then $(Y,w)$ is a relative $t$-design of $H(n,2)$ if and only if $(V,\mathcal B_{r_\nu})$ is a combinatorial $(t-1)$-$(n,r_\nu,\lambda^{(r_\nu)}_{t-1})$ design for $\nu=1,2$ and the following equality holds:
\begin{eqnarray}
&&\sum_{\nu=1}^2w_{r_\nu}\lambda^{(r_\nu)}_t(i_1,\ldots,i_t)=\sum_{\nu=1}^2N_{r_\nu}w_{r_\nu}\prod_{j=0}^{t-1}\frac{r_\nu-j}{n-j}
\label{equ:3.1}
\end{eqnarray}
for any distinct $t$ points $i_1,\ldots,i_t$ in $V$.
In above $\lambda^{(r_\nu)}_t(i_1,\ldots,i_t)$ denotes the number of blocks in $\mathcal B_{r_\nu}$ containing $\{i_1,\ldots,i_t\}$.
\label{theo:3.4}\end{theo}
\begin{theo}
Let $(V,\mathcal B_r)$ be a combinatorial $2e$-$(n,r,\lambda_{2e})$ design and $(V,\mathcal B_{n-r})$ be the complementary design of $(V,\mathcal B_r)$, with $\mathcal B_{n-r}=\{V\backslash B\mid B\in \mathcal B_r\}$.
Let $\mathcal B=\mathcal B_r\cup \mathcal B_{n-r}$ and $Y=Y_{r}\cup Y_{n-r}$ be the subset of $X$ on two shells $X_{r}\cup X_{n-r}$ corresponding to the block set $\mathcal B$.
Then $(Y,w)$ is a relative $(2e+1)$-design with constant weight, namely, $w(y)=w$ for any $y\in Y$.
\label{theo:3.5}
\end{theo}
In general, for an odd integer $t$, we do not have a natural lower bound for relative $t$-designs on $p$ shells.
 However, if $p=2$ and weight is constant on each shell, then Theorem \ref{theo:3.3} implies that $Y_{r_1}$ and $Y_{r_2}$ have the structures of combinatorial $2e$-designs.
Therefore we must have
\begin{eqnarray}
|Y|=|Y_{r_1}|+|Y_{r_2}|\geq 2{n\choose e}.
\label{equ:3.2}
\end{eqnarray}
If a relative $(2e+1)$-design $(Y,w)$ satisfies equality in \eqref{equ:3.2}, then we say $(Y,w)$ is tight.

In the case of $t=3$, Theorem \ref{theo:3.3} and \ref{theo:3.4} imply that $(Y,w)$ is a relative 3-design on two shells $X_{r_1}\cup X_{r_2}$ with constant weight on each shell if and only if the corresponding designs $(V,\mathcal B_{r_1})$ and $(V,\mathcal B_{r_2})$ are combinatorial 2-designs.
We have the following theorem.
\begin{theo}
Let $(Y,w)$ be a weighted subset of a union of two shells $X_{r_1}\cup X_{r_2}$ in $H(n,2)$.
Assume $w\equiv 1$, $n=4u-1$, $r_1=2u-1$ and $r_2=2u$.
Then $(Y,w)$ is a tight relative 3-design if and only if the corresponding design $(V,\mathcal B_{r_1})$ is a symmetric $2$-$(4u-1,2u-1,u-1)$ design and $(V,\mathcal B_{r_2})$ is the complementary design of $(V,\mathcal B_{r_1})$.
\label{theo:3.6}
\end{theo}
{\bf Proof.~}
Let $t=3$, $w_{r_1}=w_{r_2}\equiv 1$,  $n=4u-1$, $r_1=2u-1$, $r_2=2u$ in \eqref{equ:3.1}, then we have $$\lambda^{(r_1)}_3(i_1,i_2,i_3)+\lambda^{(r_2)}_3(i_1,i_2,i_3)=u-1.$$
Let $\infty$ be a point not in $V$ and $V^+=V\cup\{\infty\}$.
Let $\mathcal B^+_{r_1}=\{B\cup\{\infty\}\mid B\in \mathcal B\}$.
Then $(V,\mathcal B^+_{r_1}\cup \mathcal B_{r_2})$ is a $3$-$(4u,2u,u-1)$ Hadamard design with $8u-2$ blocks.
It is well known that the complement of any block of $3$-$(4u,2u,u-1)$ Hadamard design is again a block (cf. \cite{Driessen-1978}, Lemma 4.1).
This completes the proof.\hfill\qed\\

We give the proof of our main results, Theorem \ref{theo:3.3}, Theorem \ref{theo:3.4}, Theorem \ref{theo:3.5} in the following section. \\

For tight relative $2e$-designs in Q-polynomial association scheme $\mathfrak X$, it is known that if the stabilizer $G_0$ of $u_0$ acts transitively on each shell $X_r$, $1\leq r\leq n$, then weight function $w$ is constant on each $Y_r$, $1\leq r\leq p$ (see \cite{Bannai-B-B-2014}).
But for odd integer $t=2e+1$, we only have natural lower bound for the case $p=2$, assuming $w$ is constant on each shell.
The following problems would be interesting.
\begin{enumerate}
\item For relative $(2e+1)$-design of $H(n,2)$ on two shells, can we prove $|Y|\geq 2{n\choose e}$ holds without assuming that the weight is constant on each shell ?
\item Can we generalize Kageyama's Theorem for the relative $t$-design on $p$ shells with $p\geq 3$ ?
\item Can we generalize Kageyama's Theorem for $H(n,q)$, $q\neq 2$ ? \end{enumerate}

\section{Proof of main results}
{\bf Proof of Theorem \ref{theo:3.3}}

Let $(V,\mathcal B,w)$ be a $t$-wise balanced design.
By assumption $B\in \mathcal B$ has either $r_1$ points or $r_2$ points.
We assume $2\leq r_1<r_2\leq n-2$ to avoid the trivial cases.
Let $\mathcal B_{r_\nu}$ be the set of blocks of size $r_\nu$ for $\nu=1,2$.
Let $A$ be the incidence matrix of $\mathcal B$, i.e., $A$ is a matrix indexed by $V\times \mathcal B$ whose entry is defined by
$$A(i,B)=\left\{ \begin{array}{ll}
1 & \mbox{if $i \in B$,}\\
0 & \mbox{otherwise.}
\end{array}\right.$$
Let $A_{r_\nu}$ $(\nu=1,2)$ be a matrix indexed by $V\times \mathcal B$ defined by
$$A_{r_\nu}(i,B)=\left\{
\begin{array}{ll}
1 & \mbox{if $B\in \mathcal B_{r_\nu}$ and $i \in B$,}\\
0 & \mbox{otherwise.} \end{array}\right.$$
 By definition we have
$$A=A_{r_1}+A_{r_2}.$$
 In particular for $B\in \mathcal B_{r_\nu}~(\nu=1,2)$ we have
$$A(i,B)=A_{r_1}(i,B)+A_{r_2}(i,B)=A_{r_\nu}(i,B)$$
 and
$$\sum_{i=1}^nA(i,B)=\sum_{i=1}^n A_{r_\nu}(i,B)=r_\nu.$$
 Let $i_1,\ldots,i_{t-1}$ be distinct $t-1$ points in $V$.
 Then we have
 \begin{eqnarray}
 &&\sum_{B\in \mathcal B}w(B)\bigg\{r_2-\sum_{i=1}^n A(i,B)\bigg\}A(i_1,B)\cdots A(i_{t-1},B)\nonumber\\
 &&=\sum_{B\in \mathcal B_{r_1}}w(B)\bigg\{r_2-\sum_{i=1}^n A(i,B)\bigg\}A(i_1,B)\cdots A(i_{t-1},B)\nonumber\\
 && \quad +\sum_{B\in \mathcal B_{r_2}}w(B)\bigg\{r_2-\sum_{i=1}^n A(i,B)\bigg\}A(i_1,B)\cdots A(i_{t-1},B)\nonumber\\
 &&=(r_2-r_1)\sum_{B\in \mathcal B_{r_1}}w(B)A_{r_1}(i_1,B)\cdots A_{r_1}(i_{t-1},B)\nonumber\\
 &&=(r_2-r_1)\sum_{B\in \mathcal B_{r_1}\atop \{i_1,\ldots,i_{t-1}\} \subset B}w(B)\nonumber\\
 &&=w_{r_1}(r_2-r_1)\left|\{B\in \mathcal B_{r_1}\mid  \{i_1,\ldots,i_{t-1}\} \subset B\}\right|.
 \label{equ:4.1}\end{eqnarray}
 On the other hand we have
 \begin{eqnarray}
 && \sum_{B\in \mathcal B}w(B)\bigg\{r_2-\sum_{i=1}^n A(i,B)\bigg\}A(i_1,B)\cdots A(i_{t-1},B)\nonumber\\
 &&=\sum_{B\in \mathcal B}w(B)r_2A(i_1,B)\cdots A(i_{t-1},B)-\sum_{i=1}^n\sum_{B\in \mathcal B}w(B) A(i,B)A(i_1,B)\cdots A(i_{t-1},B)\nonumber\\
 &&=r_2\sum_{B\in \mathcal B\atop { \{i_1,\ldots,i_{t-1}\}\subset B }}w(B)-\sum_{i\not \in \{i_1,\ldots,i_{t-1}\}}\sum_{B\in \mathcal B\atop{\{i,i_1,\ldots,i_{t-1}\}\subset B}}w(B) -\sum_{i \in \{i_1,\ldots,i_{t-1}\}}\sum_{B\in \mathcal B\atop{\{i_1,\ldots,i_{t-1}\}\subset B}}w(B)\nonumber\\
 &&=r_2\lambda_{t-1}-(n-t+1)\lambda_t-(t-1)\lambda_{t-1}\nonumber\\
  &&=(r_2-t+1)\lambda_{t-1}-(n-t+1)\lambda_t.
 \label{equ:4.2}\end{eqnarray}
 Then \eqref{equ:4.1} and \eqref{equ:4.2} imply
$$
|\{B\in \mathcal B_{r_1}\mid  \{i_1,\ldots,i_{t-1}\} \subset B\}|=\frac{(r_2-t+1)\lambda_{t-1}-(n-t+1)\lambda_t}{(r_2-r_1)w_{r_1}}.
$$
 Hence $\mathcal B_{r_1}$ is a combinatorial $(t-1)$-$(n,r_1,\lambda^{(r_1)}_{t-1})$ design with $\lambda^{(r_1)}_{t-1}=\frac{(r_2+t-1)\lambda_{t-1}+(n-t+1)\lambda_t}{(r_2-r_1)w_{r_1}}$.
 Similarly we can prove the statements for $\mathcal B_{r_2}$.
 This completes the proof.\hfill\qed\\

 To give the proof of Theorem \ref{theo:3.4} and Theorem \ref{theo:3.5}, we need some preparation.
 Let $(V,\mathcal B)$ be combinatorial $(t-1)$-$(n,r,\lambda_{t-1})$ design.
 Let $i_1,\ldots,i_s$ be distinct $s$ points in $V$ and $\ell$ be an integer satisfying $0\leq \ell\leq s$.
 We define
$$ p_{\mathcal B}(\ell;i_1,\ldots,i_s)=|\{B\in \mathcal B: |B\cap \{i_1,\ldots,i_s\}| =\ell\}|.$$
 When we consider relative $t$-design, we use the following expressions.
 Let $\mathcal B_{r_\nu}$ be the block set corresponding to $Y_{r_\nu}\subset X_{r_\nu}$ (a shell in $H(n,2)$) for $\nu=1,2$.
 Let $u_1,\ldots,u_s$ be distinct $s$ points on the shell $X_1$ in $H(n,2)$.
 Let $i_j$ be the coordinate of $u_j$ which takes the value 1 for $j=1,\ldots,s$.
 Then we have the following equality
$$\big|\big\{y\in Y_{r_\nu}\mid~ \sharp\{u_i\mid (y,u_i)\in R_{r_\nu-1}\}=\ell\big\}\big|=p_{\mathcal B_{r_\nu}}(\ell; i_1,\ldots, i_s)$$
 for $\nu=1,2$.
 In the following firstly we consider some properties of $p_{\mathcal B}(\ell;i_1,\ldots,i_s)$ for $(t-1)$-$(n,r,\lambda_{t-1})$ designs.
 It is well known that, for a combinatorial $(t-1)$-$(n,r,\lambda_{t-1})$ design $\mathcal B$, the following equality holds:
 $$p_{\mathcal B}(s;i_1,\ldots,i_s)=\lambda_s=\frac{{n-s\choose r-s}}{{n\choose r}}|\mathcal B|$$
 for $s=0,\ldots,t-1$.
 We also use the notation $\lambda^{\mathcal B}_t(i_1,\ldots,i_t)$ instead of $p_{\mathcal B}(t;i_1,\ldots,i_t)$ for a combinatorial $(t-1)$-design $\mathcal B$.
 Before beginning the arguments on the designs we firstly introduce a combinatorial formula which we use to prove our theorem, although it may be well known already.
\begin{lem} \label{lem:4.1}
 For non-negative integers $\alpha, \beta$ and $\gamma$ satisfying $\alpha\leq \beta$, the following equality holds.
 \begin{eqnarray}
 &&\sum_{j=0}^{\alpha}(-1)^j{\alpha\choose j}{\beta+\gamma-j\choose \beta}={\beta-\alpha+\gamma\choose \gamma}.
\label{equ:4.3}
\end{eqnarray}
\end{lem}
{\bf Proof. ~}
We have the following formulas for polynomials of $x$.
$$(1-x)^\alpha=\sum_{j=0}^\alpha (-1)^j{\alpha\choose j}x^j,$$
$$(1-x)^{-\beta-1}=\sum_{i=0}^{\infty}{\beta+i\choose i}x^i.$$
Hence we have
$$(1-x)^{\alpha-\beta-1}=\sum_{i=0}^{\infty}\sum_{j=0}^\alpha (-1)^j{\alpha\choose j}{\beta+i\choose\beta}x^{i+j}.$$
Let $i+j=\gamma$. Then we have
\begin{eqnarray}
&&(1-x)^{\alpha-\beta-1}=\sum_{\gamma=0}^{\infty}\sum_{j=0}^\alpha (-1)^j{\alpha\choose j}{\beta+\gamma-j\choose\beta}x^{\gamma}.
\label{equ:4.4}
\end{eqnarray}
On the both hand since $\beta-\alpha\geq 0$, we have
\begin{eqnarray}
&&(1-x)^{\alpha-\beta-1}=(1-x)^{-(\beta-\alpha)-1}=\sum_{\gamma=0}^{\infty}{\beta-\alpha+\gamma \choose \gamma}x^\gamma.
\label{equ:4.5}
\end{eqnarray}
\eqref{equ:4.4} and \eqref{equ:4.5} complete the proof.\hfill\qed\\

%Proposition 4.2
\begin{pro} Definition and notation are given as above.
Let $(V,\mathcal B)$ be a combinatorial $(t-1)$-$(n,r,\lambda_{t-1})$ design and $N=|\mathcal B|$.
\begin{enumerate}
\item Let $\{i_1,\ldots,i_s\}$ be an $s$-element subset of $V$.
Then the following equation holds for $s$ and $\ell$ satisfying $1\leq s\leq  t-1$ and $0\leq \ell\leq s$.
\begin{eqnarray}
p_{\mathcal B}(\ell;i_1,\ldots,i_s)=\frac{{n-s \choose r-\ell}}{{n\choose r}}N.
\label{equ:4.6}
\end{eqnarray}
\item Let $\{i_1,\ldots,i_t\}$ be a $t$-element subset of $V$.
Then the following equality holds:
\begin{equation}
p_\mathcal B(\ell;i_1,\ldots,i_t)
=\frac{N}{{n\choose r}}\left\{{n-t\choose r-\ell}
-(-1)^{t-\ell}{n-t \choose r-t}  \right\}
+(-1)^{t-\ell} \lambda^{\mathcal B}_t(i_1,\ldots,i_t)
\label{equ:4.7}
\end{equation}
for any integer satisfying $0\leq \ell\leq t-1$.
%Here we denote $p_\mathcal B(t;i_1,\ldots,i_t)$ by $\lambda^{\mathcal B}_t(i_1,\ldots,i_t)$ for simplicity.
\end{enumerate}
\label{pro:4.2}\end{pro}
{\bf Proof. ~}
The formula \eqref{equ:4.6} is well known (cf. \cite{Driessen-1978}).
The formula \eqref{equ:4.7} also might be already well known, however we will give the proof below.
By inclusion-exclusion principle, we have
\begin{eqnarray}
&&p_{\mathcal B}(\ell;i_1,\ldots,i_t)
=\lambda_{\ell}-{t-\ell\choose 1}\lambda_{\ell+1}
+{t-\ell\choose 2}\lambda_{\ell+2}+\cdots
+(-1)^{t-\ell-1}{t-\ell\choose t-\ell-1}\lambda_{t-1}\nonumber\\
&& \quad +(-1)^{t-\ell}\lambda^{\mathcal B}_t(i_1,\ldots,i_t)\nonumber\\
&&=\sum_{j=0}^{t-\ell-1}(-1)^j{t-\ell\choose j}\lambda_{\ell+j}+(-1)^{t-\ell}\lambda^{\mathcal B}_t(i_1,\ldots,i_t).
\label{equ:4.8}\end{eqnarray}
Since $(V,\mathcal B)$ is a combinatorial $(t-1)$-$(n,r,\lambda_{t-1})$ design, the following equation is known.
\begin{eqnarray}
&&\lambda_{\ell+j}=\frac{N}{{n\choose r}}{n-(\ell+j)\choose n-r}\quad \mbox{for $\ell+j\leq t-1 $}.
\label{equ:4.9}\end{eqnarray}
Then \eqref{equ:4.8} and \eqref{equ:4.9} imply
\begin{eqnarray}
&&p_{\mathcal B}(\ell;i_1,\ldots,i_t)=\frac{N}{{n\choose r}}\sum_{j=0}^{t-\ell-1}(-1)^j{t-\ell\choose j}{n-\ell-j\choose n-r}+(-1)^{t-\ell}\lambda^{\mathcal B}_t(i_1,\ldots,i_t)
\nonumber\\
&&=\frac{N}{{n\choose r}}\sum_{j=0}^{t-\ell}(-1)^j{t-\ell\choose j}{n-\ell-j\choose n-r}-(-1)^{t-\ell}\frac{N}{{n\choose r}}{n-t\choose n-r}+(-1)^{t-\ell}\lambda^{\mathcal B}_t(i_1,\ldots,i_t)\nonumber\\
&&=\frac{N}{{n\choose r}}\sum_{j=0}^{t-\ell}(-1)^j{t-\ell\choose j}{(n-r)+(r-\ell) -j\choose n-r}\nonumber\\
&&\quad -(-1)^{t-\ell}\frac{N}{{n\choose r}}{n-t\choose n-r}+(-1)^{t-\ell}\lambda^{\mathcal B}_t(i_1,\ldots,i_t).
\label{equ:4.10}\end{eqnarray}
Then apply formula \eqref{equ:4.3} with $\alpha=t-\ell$, $\beta=n-r$ and $\gamma=r-\ell$, to equation \eqref{equ:4.10}
we have
$$
p_{\mathcal B}(\ell;i_1,\ldots,i_t)
=\frac{N}{{n\choose r}}\bigg\{{n-t\choose r-\ell}
-(-1)^{t-\ell}{n-t\choose n-r}\bigg\}+(-1)^{t-\ell}\lambda^{\mathcal B}_t(i_1,\ldots,i_t).
$$
This completes the proof of Proposition \ref{pro:4.2} (2).\hfill\qed\\

Let $(V,\mathcal B)$ be a combinatorial $(t-1)$-$(n,r, \lambda_{t-1})$ design.
Let $\mathcal B^c=\{V\backslash B\mid B\in \mathcal B\}$.
Then it is known that $(V,\mathcal B^c)$ is a $(t-1)$-$(n,n-r,\lambda^c_{t-1})$.
\begin{cor}
Notations and definitions are given as above.
The following equality holds
\begin{eqnarray}
\lambda^{\mathcal B^c}_t(i_1,\ldots,i_t)-(-1)^t \lambda^{\mathcal B}_t(i_1,\ldots,i_t)=\frac{N}{{n\choose r}}\bigg\{{n-t\choose r}-(-1)^{t}{n-t\choose n-r}\bigg\}.
\label{equ:4.11}\end{eqnarray}
\label{cor:4.3}
\end{cor}
{\bf Proof. ~}
Proposition \ref{pro:4.2} (2) and $\lambda^{\mathcal B^c}_t(i_1,\ldots,i_t)=p_{\mathcal B}(0;i_1,\ldots,i_t)$ implies \eqref{equ:4.11}.
\hfill\qed\\

Now we are ready to discuss relative $t$-design $(Y,w)$ in $H(n,2)$ on two shells $X_{r_1}\cup X_{r_2}$ and give the proof of Theorem \ref{theo:3.4} we stated in Section \ref{sec:main}.
Let $Y_{r_\nu}=Y\cap X_{r_\nu}$ for $\nu=1,2$.
As we have seen in Section \ref{sec:main}, $Y_{r_1}$ and $Y_{r_2}$ have the structure of combinatorial $(t-1)$-$(n,r_\nu,\lambda^{(r_\nu)}_{t-1})$ design, i.e.,
$(V,\mathcal B_{r_\nu})$ is a combinatorial $(t-1)$-$(n,r_\nu,\lambda^{(r_\nu)}_{t-1})$ design for $\nu=1,2$.
Where $\mathcal B_{r_\nu}=\{B_y\mid y\in Y_{r_\nu}\}$, $B_y=\{i\mid 1\leq i\leq n,~ y_i=1\}$.\\

\noindent{\bf Proof of Theorem \ref{theo:3.4}}

Let $u_1,\ldots,u_s$ be distinct $s$ points in $X_1$.
Let us consider the defining equation of relative $t$-designs.
\begin{eqnarray}
\sum_{\nu=1}^2\frac{N_{r_\nu}w_{r_\nu}}{{n\choose r_\nu}}\sum_{x\in X_{r_\nu}}\prod_{j=1}^s\phi^{(1)}_{u_j}(x)
=\sum_{\nu=1}^2\sum_{y\in Y_{r_\nu}}w_{r_\nu}\prod_{j=1}^s\phi^{(1)}_{u_j}(y).
\label{equ:4.12}\end{eqnarray}
Let $\left(Q_j(\ell)\right)_{0\leq \ell,j\leq n}$ be the second eigenmatrix of $H(n,2)$.
It is well known that $Q_1(\ell)=n-2\ell$ holds for $0\leq \ell\leq n$.
Since $(x,u_j)\in R_{r_\nu-1}\cup R_{r_\nu+1}$ holds for any $x\in X_{r_\nu}$, we have the following equality on the cardinality of the set for each $\ell$
satisfying $0\leq \ell \leq s$.
$$
\sharp\bigg\{x\in X_{r_\nu}\bigg| \begin{array}{l}|\{j\mid (x,u_j)\in R_{r_\nu-1}\}|=\ell, ~\mbox{and}\\
|\{j\mid(x,u_j)\in R_{r_i+1}\}|
=s-\ell\end{array}\bigg\}={s\choose \ell}{n-s\choose r_\nu-\ell}.
$$
Therefore the left hand side of \eqref{equ:4.12} equals
\begin{eqnarray}
&&\sum_{\nu=1}^2\frac{N_{r_\nu}w_{r_\nu}}{{n\choose r_\nu}}\sum_{\ell=0}^s{s\choose \ell}{n-s\choose r_\nu-\ell}Q_1(r_\nu-1)^{\ell}Q_1(r_\nu+1)^{s-\ell}.
\label{equ:4.13}\end{eqnarray}
Next we consider the right hand side of \eqref{equ:4.12}.
Let $i_j$ be the coordinate of $u_i$ whose entry is $1$ for $1\leq i\leq t$.
Let $\mathcal B_{r_\nu}$ be the block set corresponding to $Y_{r_\nu}$.
By \eqref{equ:4.2} and \eqref{equ:4.6}, if $s<t$, then the right hand side of  \eqref{equ:4.12} equals to the following formula.
 \begin{eqnarray}
 && \sum_{i=1}^2w_{r_\nu}\sum_{\ell=0}^s{s\choose \ell} p_{\mathcal B_{r_\nu}}(\ell;i_1,\ldots,i_s) Q_1(r_\nu-1)^{\ell}Q_1(r_\nu+1)^{s-\ell}\nonumber\\
&&= \sum_{i=1}^2w_{r_\nu}\sum_{\ell=0}^s{s\choose \ell} \frac{{n-s\choose r_\nu-\ell}}{{n\choose r_\nu}} Q_1(r_\nu-1)^{\ell}Q_1(r_\nu+1)^{s-\ell}. \nonumber
\end{eqnarray}
Thus the equality \eqref{equ:4.12} holds for any $s=1,\ldots,t-1$.
Next, let $s=t$.
 Then Proposition \ref{pro:4.2} (2) implies that the right hand side of \eqref{equ:4.12} equals to the following formula
\begin{eqnarray}
&& \sum_{\nu=1}^2w_{r_\nu}\sum_{\ell=0}^t{t\choose \ell} \bigg(\frac{N_{r_\nu}}{{n\choose r_\nu}}\left\{{n-t\choose r_\nu-\ell}-(-1)^{t-\ell}{n-t \choose r_\nu-t}  \right\}+(-1)^{t-\ell} \lambda^{(r_\nu)}_t(i_1,\ldots,i_t)\bigg)\times
\nonumber\\
&&Q_1(r_\nu-1)^{\ell}Q_1(r_\nu+1)^{t-\ell}.
\label{equ:4.14}\end{eqnarray}
Here we use the notation $\lambda^{(r_\nu)}_t(i_1,\ldots,i_t)$ instead of $\lambda^{B_{r_\nu}}_t(i_1,\ldots,i_t)~(=p_{\mathcal B_{r_\nu}}(t;i_1,\ldots,i_t))$
for simplicity.
Then \eqref{equ:4.13} and \eqref{equ:4.14} imply that for $s=t$, \eqref{equ:4.12} is equivalent to the following equation.
\begin{eqnarray}
&&\sum_{\nu=1}^2\frac{N_{r_\nu}w_{r_\nu}}{{n\choose r_\nu}}\sum_{\ell=0}^t{t\choose \ell}{n-t\choose r_\nu-\ell}Q_1(r_\nu-1)^{\ell}Q_1(r_\nu+1)^{t-\ell}
\nonumber\\
&&= \sum_{\nu=1}^2w_{r_\nu}\sum_{\ell=0}^t{t\choose \ell} \bigg(\frac{N_{r_\nu}}{{n\choose r_\nu}}\left\{{n-t\choose r_\nu-\ell}-(-1)^{t-\ell}{n-t \choose r_\nu-t}  \right\}+(-1)^{t-\ell} \lambda^{(r_\nu)}_t(i_1,\ldots,i_t)\bigg)\times \nonumber\\
&&\quad Q_1(r_\nu-1)^{\ell}Q_1(r_\nu+1)^{t-\ell} \nonumber\\
&&= \sum_{\nu=1}^2 \frac{N_{r_\nu}w_{r_\nu}}{{n\choose r_\nu}}\sum_{\ell=0}^t{t\choose \ell}{n-t\choose r_\nu-\ell}Q_1(r_\nu-1)^{\ell}Q_1(r_\nu+1)^{t-\ell}
\nonumber\\
&&\quad - \sum_{\nu=1}^2 \frac{N_{r_\nu}w_{r_\nu}}{{n\choose r_\nu}}{n-t \choose r_\nu-t} \sum_{\ell=0}^t(-1)^{t-\ell}{t\choose \ell}Q_1(r_\nu-1)^{\ell}Q_1(r_\nu+1)^{t-\ell} \nonumber\\
&&\quad + \sum_{\nu=1}^2w_{r_\nu}p_{\mathcal B_{r_\nu}}(t;i_1,\ldots,i_t)\sum_{\ell=0}^t(-1)^{t-\ell}{t\choose \ell}Q_1(r_\nu-1)^{\ell}Q_1(r_\nu+1)^{t-\ell}
\nonumber\\
&&= \sum_{\nu=1}^2 \frac{N_{r_\nu}w_{r_\nu}}{{n\choose r_\nu}}\sum_{\ell=0}^t{t\choose \ell}{n-t\choose r_\nu-\ell}Q_1(r_\nu-1)^{\ell}Q_1(r_\nu+1)^{t-\ell}
\nonumber\\
&&\quad - \sum_{\nu=1}^2 \frac{N_{r_\nu}w_{r_\nu}}{{n\choose r_\nu}}{n-t \choose r_\nu-t}\big(Q_1(r_\nu-1)-Q_1(r_\nu+1)\big)^t\nonumber\\
&&\quad+ \sum_{\nu=1}^2w_{r_\nu}\lambda^{(r_\nu)}_t(i_1,\ldots,i_t)\big(Q_1(r_\nu-1)-Q_1(r_\nu+1)\big)^t.
\label{equ:4.15}\end{eqnarray}
Since $Q_1(r_\nu-1)-Q_1(r_\nu+1)=4$, \eqref{equ:4.15} implies
\begin{eqnarray}
&& \sum_{\nu=1}^2w_{r_\nu} \lambda^{(r_\nu)}_t(i_1,\ldots,i_t)=\sum_{\nu=1}^2\frac{N_{r_\nu}w_{r_\nu}}{{n\choose r_\nu}}{n-t\choose r_\nu-t}=\sum_{\nu=1}^2N_{r_\nu}w_{r_\nu}\prod_{j=0}^{t-1}\frac{r_\nu-j}{n-j}.
\label{equ:4.16}\end{eqnarray}
Thus we proved that equation \eqref{equ:4.12} with $s=t$ is equivalent to \eqref{equ:3.1}.
This completes the proof.\hfill \qed\\

\noindent{\bf Proof of Theorem \ref{theo:3.5}}

Let $r_1=r$ and $r_2=n-r$.
Since $(V, \mathcal B_{r_1})$ and $(V, \mathcal B_{r_2})$ are combinatorial $2e$-$(n,r_\nu,\lambda^{(r_\nu)}_{2e})$ designs which are complementary designs of each other $N_{r_1}=N_{r_2}$ holds and Theorem \ref{theo:3.4} implies that it is enough that we prove the equation \eqref{equ:3.1} holds for $t=2e+1$.
In the proof of Theorem \ref{theo:3.4} it is shown that \eqref{equ:3.1} is equivalent to \eqref{equ:4.16}.
On the other hand Corollary \ref{cor:4.3} implies
\begin{eqnarray}
&&\lambda^{(r_1)}_{2e+1}(i_1,\ldots,i_{2e+1})+\lambda^{(r_2)}_{2e+1}(i_1,\ldots,i_{2e+1})\nonumber\\
&&=\frac{N}{{n\choose r}}\bigg({n-2e-1\choose r}-(-1)^{2e+1}{n-2e-1\choose n-r}\bigg)\nonumber\\
&&=\frac{N}{{n\choose r}}\bigg({n-2e-1\choose r}+{n-2e-1\choose n-r}\bigg).
\label{equ:4.17}\end{eqnarray}
If $t=2e+1$, $r_1=r$ and $r_2=n-r$, then \eqref{equ:4.16} is equivalent to \eqref{equ:4.17}.
Moreover, this implies  $w_{r_1}=w_{r_2}$.
\hfill\qed\\

\noindent
The following proposition is very useful.
\begin{pro}
Let $(Y,w)$ be a tight relative $3$-design on two shells $X_{r_1}\cup X_{r_2}$ in $H(n,2)$. Assume $r_1+r_2=n,$ $r_1<r_2$ and $w_{r_1}=w_{r_2}$.
Then
$$\lambda^{(r_1)}_3(i_1,i_2,i_3)\geq 1$$
for any 3-point subset $\{i_1,i_2,i_3\}\subset V\backslash B$, with any $B\in \mathcal B_{r_2}$.
\label{pro:4.4}
 \end{pro}

\noindent{\bf Proof. ~}
By assumption $w_{r_1}=w_{r_2}$, $r_1+r_2=n$ and by Theorem \ref{theo:3.2} and Theorem \ref{theo:3.3} the combinatorial designs $(V, \mathcal B_{r_1})$ and $(V, \mathcal B_{r_2})$ corresponding to $Y_{r_1}$ and $Y_{r_2}$ are symmetric $2$-$(n,r_1,\lambda^{(r_1)}_{2})$ and $2$-$(n,r_2,\lambda^{(r_2)}_{2})$ designs respectively.
Therefore $|\mathcal B_{r_1}|=|\mathcal B_{r_2}|=n$ and \eqref{equ:3.1} imply
\begin{eqnarray}
&&\lambda^{(r_1)}_3(i_1,i_2,i_3)+\lambda^{(r_2)}_3(i_1,i_2,i_3) \nonumber\\
&&=\frac{1}{(n-1)(n-2)}\bigg(r_1(r_1-1)(r_1-2)+r_2(r_2-1)(r_2-2)\bigg)\nonumber\\
&&=\frac{n^2-3nr_2+3r_2^2-n}{n-1} \nonumber
\end{eqnarray}
for any 3-point subset $\{i_1,i_2,i_3\}\subset V$.
Let $B=\{a_1,\ldots,a_{r_2}\}\in   \mathcal B_{r_2}$.
Let $V\backslash B=\{a_{r_2+1},\ldots,a_n\}$.
Assume that $\lambda^{(r_1)}_3(a_{r_2+1},a_{r_2+2},a_{r_2+3})=0$ for $\{a_{r_2+1},a_{r_2+2},a_{r_2+3}\}\subset V\backslash B$.
Then
$$\lambda^{(r_2)}_3(a_{r_2+1},a_{r_2+2},a_{r_2+3})=\frac{n^2-3nr_2+3r_2^2-n}{n-1}$$
holds.
Let
\begin{eqnarray}\alpha_3=\frac{n^2-3nr_2+3r_2^2-n}{n(n-1)}.
\label{equ:4.18}
\end{eqnarray}
Count the number of blocks in $\mathcal B_{r_2}$ according to the manner given below.
Note that the following formula for symmetric design are well known
\begin{eqnarray}
&&\lambda^{(r_2)}_2=\frac{r_2(r_2-1)}{n-1},\quad \lambda^{(r_2)}_1=r_2.\label{equ:4.19}
\end{eqnarray}
Then we have the following equations.
\begin{eqnarray}&&
\big|\big\{B\in \mathcal B_{r_2}\big|
~|B\cap \{a_{r_2+1},a_{r_2+2},a_{r_2+3}\}|=3\big\}\big|
=\alpha_3,\label{equ:4.20}\\
&&\big|\big\{B\in \mathcal B_{r_2}\big |
|\{B\cap\{a_{r_2+1},a_{r_2+2}, a_{r_2+3}\}|=2\big\}\big|
=\lambda^{(r_2)}_2-\alpha_3,\label{equ:4.21}\\
&&\big|\big\{B\in \mathcal B_{r_2}\big |
|\{B\cap \{a_{r_2+1},a_{r_2+2}, a_{r_2+3}\}|=1\big\}\big|
=r_2-\lambda^{(r_2)}_2.
\label{equ:4.22}
\end{eqnarray}
By \eqref{equ:4.18}--\eqref{equ:4.22}, we have
\begin{eqnarray}
&&n=|\mathcal B_{r_2}|\geq \alpha_3+3(\lambda^{(r_2)}_2-\alpha_3)+
3(r_2-\lambda^{(r_2)}_2)
\nonumber\\
&&=3r_2-2\frac{(n^2-3nr_2+3r_2^2-n)}
{n-1}
\nonumber\\
&&=\frac{-2n^2+2n+9nr_2-3r_2-6r_2^2}{n-1}.
\label{equ:4.23}
\end{eqnarray}
\eqref{equ:4.23} implies
$$\frac{3(n-r_2)(n-2r_2-1)}{n-1}\geq 0.$$
Hence we must have $2r_2\leq n-1$.
On the other hand by the assumption, we have $n-r_2=r_1<r_2$, this is a contradiction.
This completes the proof of Proposition {\ref{pro:4.4}.\hfill\qed

\section{Tight relative $(2e+1)$-designs}\label{sec:design_odd}
Let $(Y,w)$ be a tight relative $(2e+1)$-design on $X_{r_1}\cup X_{r_2}$.
We assume weight $w$ is constant on each $Y_{r_\nu}$ and let $w_{r_\nu}=w(y)$ for $y\in Y_{r_\nu}$, $\nu=1,2$.
Then $Y_{r_\nu}=Y\cap X_{r_\nu}$ is a tight $2e$-design for $\nu=1,2$.
Namely, $N_{r_\nu}=|Y_{r_\nu}|={n\choose e}$.
If $r_1=1$ or $r_2=n-1$, then we must have trivial case $Y_{r_1}=X_{r_1}$ or $Y_{r_2}=X_{r_2}$.
Hence in the following we assume $2\leq r_1<r_2\leq n-2$.
%%%%%%%%%%%

\subsection{Tight relative $3$-designs}

\noindent
{\bf Method to get feasible parameters}\\
The formula given in Theorem \ref{pro:4.2} for $t=3$ give the following formulas.
In the following we use the notation $\lambda^{(r_\nu)}_3(i_1,i_2,i_3)$ instead of $p_{\mathcal B_{r_\nu}}(3;i_1,i_2,i_3)$ for simplicity.
\begin{eqnarray}
&&p_{\mathcal B_{r_\nu}}(0;i_1)=\lambda^{(r_\nu)}_0 =n-r_\nu,
\quad p_{\mathcal B_{r_\nu}}(1;i_1)=\lambda^{(r_\nu)}_1  =\frac{n{n-1\choose r_\nu-1}}{{n\choose r_\nu}}=r_\nu,\nonumber\\
&&p_{\mathcal B_{r_\nu}}(0;i_1,i_2)=n-r_\nu-\frac{r_\nu(n-r_\nu)}{n-1},
\quad p_{\mathcal B_{r_\nu}}(1;i_1,i_2)=\frac{r_\nu(n-r_\nu)}{n-1},\nonumber\\
&&p_{\mathcal B_{r_\nu}}(2;i_1,i_2)=\lambda^{(r_\nu)}_2=\frac{(r_\nu-1)r_\nu}{n-1}.\label{equ:5.1}
\end{eqnarray}
For any distinct 3 points $i_1,i_2,i_3\in V$, the following equality holds.
\begin{eqnarray}
&&\sum_{\nu=1}^2w_{r_\nu}\lambda^{(r_\nu)}_3(i_1,i_2,i_3)
=\sum_{\nu=1}^2 w_{r_\nu}\frac{r_\nu(r_\nu-1)(r_\nu-2)}{(n-1)(n-2)}.
\label{equ:5.2}
\end{eqnarray}
The equations \eqref{equ:5.1} and \eqref{equ:5.2} give the equivalent condition for $(V,\mathcal B_{r_1})$ and $(V,\mathcal B_{r_2})$ to give a tight relative 3-design $Y=Y_{r_1}\cup Y_{r_2}$.
Since $(V,\mathcal B_{r_1})$ and $(V,\mathcal B_{r_2})$ cannot be combinatorial 3-designs, there exist $\{a_1,a_2,a_3\}$ and $\{b_1,b_2,b_3\}$ in $V$ satisfying
$$\lambda^{(r_1)}_2\geq \lambda^{(r_1)}_3(a_1,a_2,a_3)>\lambda^{(r_1)}_3(b_1,b_2,b_3)\geq 0.$$
Then \eqref{equ:5.2} implies
$$
 w_{r_1}(\lambda^{(r_1)}_3(a_1,a_2,a_3)-\lambda^{(r_1)}_3(b_1,b_2,b_3))=w_{r_2}(\lambda^{(r_2)}_3(b_1,b_2,b_3)-\lambda^{(r_2)}_3(a_1,a_2,a_3)).
$$
Therefore we must have
$$\lambda^{(r_2)}_2\geq \lambda^{(r_2)}_3(b_1,b_2,b_3)>\lambda^{(r_2)}_3(a_1,a_2,a_3)\geq 0.$$
Hence
$$w_{r_2}=\frac{\lambda^{(r_1)}_3(a_1,a_2,a_3)-\lambda^{(r_1)}_3(b_1,b_2,b_3)}{\lambda^{(r_2)}_3(b_1,b_2,b_3)-\lambda^{(r_2)}_3(a_1,a_2,a_3)}w_{r_1}.$$
Let
$$
\alpha=\frac{\lambda^{(r_1)}_3(a_1,a_2,a_3)-\lambda^{(r_1)}_3(b_1,b_2,b_3)}{\lambda^{(r_2)}_3(b_1,b_2,b_3)-\lambda^{(r_2)}_3(a_1,a_2,a_3)}.
$$
Then by definition $\alpha>0$, $w_{r_2}=\alpha w_{r_1}$ and
$$\begin{array}{l} \medskip
1\leq \lambda^{(r_1)}_3(a_1,a_2,a_3)-\lambda^{(r_1)}_3(b_1,b_2,b_3)\leq \lambda^{(r_1)}_2,\\
1\leq \lambda^{(r_2)}_3(b_1,b_2,b_3)-\lambda^{(r_2)}_3(a_1,a_2,a_3)\leq \lambda^{(r_2)}_2.
\end{array} $$
 The equation \eqref{equ:5.2} implies that the following holds for any distinct  $i_1,i_2,i_3$ in $V$.
$$
\lambda^{(r_2)}_3(i_1,i_2,i_3)= \frac{r_1(r_1-1)(r_1-2)}{\alpha(n-1)(n-2)}+\frac{r_2(r_2-1)(r_2-2)}{(n-1)(n-2)}-\frac{1}{\alpha}\lambda^{(r_1)}_3(i_1,i_2,i_3).
$$
 We explained how to list feasible parameters $n,r_1,r_2,N_{r_1},N_{r_2}$ and $\frac{w_{r_2}}{w_{r_1}}$ for the case $t=3$.
 For the case $t\geq 4$ we use the same method.
 As we have seen in Theorem \ref{theo:3.3}, the existence of relative $t$-design on two shells $X_{r_1}\cup X_{r_2}$ is equivalent to the existence of combinatorial $(t-1)$-$(n,r_1,\lambda^{(r_1)}_{t-1})$ and $(t-1)$-$(n,r_2,\lambda^{(r_2)}_{t-1})$ design with $\lambda^{(r_1)}_t(i_1,\ldots,i_t)$ and
 $\lambda^{(r_2)}_t(i_1,\ldots,i_t)$ satisfying the equality \eqref{equ:3.1}.
 In the following sections we often use the terminology, $\lambda_t$-sequence of a $(t-1)$-$(v,k,\lambda_{t-1})$ design $(V,\mathcal B)$.
 The definition is given as follows.
 \begin{dfe}
 Let $(V,\mathcal B)$ be a $(t-1)$-$(v,k,\lambda_{t-1})$ design.
 Let $\ell_1,\ell_2,\ldots,\ell_j$ be integers satisfying $0\leq \ell_1<\ell_2<\cdots <\ell_j\leq \lambda_{t-1}$ and
 $$\big|\big\{\{i_1,\ldots,i_t\}\subset V \big|~\mathcal \lambda^{\mathcal B}_t (i_1,i_2,\ldots,i_t)=\ell_s\big\}\big|=a_{\ell_s}>0$$
 for $s=1,\ldots,j$.
  We call the sequence $(a_{\ell_1}\ast \ell_1,\ldots,a_{\ell_j}\ast \ell_j)$ the $\lambda_t$-sequence of $(t-1)$-$(v,k,\lambda_{t-1})$ design.
 We call $j$ the length of the $\lambda_t$-sequence.
 \end{dfe}

 In the following we give the list of all the feasible parameters satisfying the integral conditions in \eqref{equ:5.1} and \eqref{equ:5.2} up to $n=100$.
  We divide the list into four cases according to the following conditions on $r_1,r_2$ and $w_{r_1},w_{r_2}$.

 {Case 1:} $r_1+r_2=n$ and $w_{r_1}=w_{r_2}$.

  {Case 2:} $r_1+r_2=n$ and $w_{r_1}\neq w_{r_2}$.

  {Case 3:} $r_1+r_2\not=n$ and $w_{r_1}= w_{r_2}$.

  {Case 4:} $r_1+r_2\not=n$ and $w_{r_1}\neq w_{r_2}$.\\

 \noindent{{\bf Case 1:} $r_1+r_2=n$ and $w_{r_1}=w_{r_2}$.}
 {\small
 \begin{center}
 	\begin{tabular}{ccc}
 $\begin{array}{|c|c|c|c|} \hline
 n & r_1 & \lambda_2^{(r_1)}&\\ \hline
 7^\ast & 3 & 1  & 1 \\  \hline
 11^\ast & 5 & 2 &  1 \\  \hline
 13 & 4 & 1  &  1\\  \hline% PG(2,3)
 15^\ast & 7 & 3 &  5\\  \hline
 16 & 6 & 2 &  3 \\  \hline
 19^\ast & 9 & 4 &  6 \\ \hline
 21 & 5 & 1 &  1 \\  \hline% PG(2,4)
 22 & 7 & 2 & \times\\  \hline
 23^\ast & 11 & 5 &  1106  \\ \hline
 25 & 9  & 3 &  78\\  \hline
 27^\ast & 13 & 6 &  208310  \\ \hline
 29 & 8  & 2 &  \times \\  \hline
 31 & 6  & 1 &  1 \\ \hline
 31 & 10 & 3 &  151 \\ \hline
 31^\ast & 15 & 7 &  10374196953\\ \hline
 34 & 12 & 4 & \times \\  \hline
 35^\ast & 17 & 8 & \geq 108131\\  \hline
 36 & 15 & 6 &  \geq 25634\\  \hline
 37 & 9 & 2  &  4\\  \hline
 39^\ast & 19 & 9 &  \geq 5.87\cdot 10^{14} \\  \hline
 40 & 13 & 4 &  \geq 1108800\\  \hline
 41 & 16 & 6  &  \geq 115307\\  \hline%  2^{10}%, \cite{Spence-1995}, \cite{Spence-1993}  \equiv -1(\mod 41)
 43 & 7 & 1 &  \times \\ \hline
 43 & 15 & 5 &  \times \\ \hline
 43^\ast & 21 & 10 &  \geq 82\\  \hline
 45 & 12 & 3 &  \geq 3752\\  \hline
 \end{array}$
  &
 $\begin{array}{|c|c|c|c|}\hline
 n & r_1 & \lambda_2^{(r_1)}&\\ \hline
 46 & 10 & 2 & \times\\  \hline
 47^\ast & 23 & 11 &  \geq 55\\  \hline
 49 & 16 & 5  &   \geq 12146 \\  \hline
 51^\ast & 25 & 12 & \geq 1\\ \hline
 52 & 18 & 6 &  \times\\ \hline
 53 & 13 & 3 &  \times\\ \hline
 55^\ast & 27 & 13 &  \geq 1\\ \hline
 56 & 11 & 2 & \geq 5\\ \hline
 57 & 8 & 1 & 1\\ \hline
 58 & 19 & 6 &  \times\\ \hline
 59^\ast & 29 & 14 &  \geq 1\\ \hline
 61 & 16 & 4 & \geq 6\\ \hline
 61 & 21 & 7 & \times\\ \hline
 61 & 25 & 10 & \geq 24\\ \hline
 63^\ast & 31 & 15 & \geq 10^{17}\\ \hline
 64 & 28 & 12 & \geq 8784\\ \hline
 66 & 26 & 10 & \geq 588\\ \hline
 67 & 12 & 2 & \times\\ \hline
 67 & 22 & 7 & \times\\ \hline
 67^\ast & 33 & 16 & \geq 1\\ \hline
 69 & 17 & 4 & \geq 4\\ \hline
 70 & 24 & 8 & \geq 28\\ \hline
 71 & 15 & 3 & \geq 72\\ \hline
 71 & 21 & 6 & \geq 2\\ \hline
 71^\ast & 35 & 17 & \geq 9\\ \hline
 73 & 9 & 1 & 1\\ \hline
 \end{array}$
 &
 $\begin{array}{|c|c|c|c|}\hline
 n & r_1 & \lambda_2^{(r_1)}&\\ \hline
 75^\ast & 37 & 18 &  \geq 1\\ \hline
 76 & 25 & 8 &  \times\\ \hline
 77 & 20 & 5 &  \times\\ \hline
 78 & 22 & 6 & \geq 3\\ \hline
 79 & 13 & 2 & \geq 2\\ \hline
 79 & 27 & 9 & \geq 1463\\ \hline
 79^\ast & 39 & 19 & \geq 2091\\ \hline
 81 & 16 & 3 & ?\\ \hline
 83^\ast & 41 & 20 & \geq 1\\ \hline
 85 & 21 & 5 & \geq 213964\\ \hline
 85 & 28 & 9 & ?\\ \hline
 85 & 36 & 15 & ?\\ \hline
 86 & 35 & 14 &  \times\\ \hline
 87^\ast & 43 & 21 &  \geq 1\\ \hline
 88 & 30 & 10 & \times\\ \hline
 89 & 33 & 12 &  \times\\ \hline
 91 & 10 & 1 & 4\\ \hline
 91 & 36 & 14 &  \times\\ \hline
 91^\ast & 45 & 22 & \geq 1\\ \hline
 92 & 14 & 2 &  \times\\ \hline
 93 & 24 & 6 &  \times\\ \hline
 94 & 31 & 10 &  \times\\ \hline
 95^\ast & 47 & 23 &  \geq 1\\ \hline
 96 & 20 & 4 & \geq 2\\ \hline
 97 & 33 & 11 & ?\\ \hline
 99^\ast & 49 & 24 & \geq 1\\ \hline
 100 & 45 & 20 &  \geq 1\\ \hline
 \end{array}$\\
 &	{Table 1} &
 \end{tabular}
 \end{center}
}
{{\bf Remark}}
\begin{enumerate}
\item  In the table given above, if a symmetric design of corresponding parameters $(n,r_1,\lambda^{(r_1)}_2)$ exists then a tight relative $3$-design on two shells $X_{r_1}\cup X_{r_2}$, with $r_2=n-r_1$, exists.
 \item  $n^\ast$ denotes the case which is 2-$(4u-1,2u-1,u-1)$ Hadamard design with $n=4u-1$.
 In this case corresponding design $(V, \mathcal B_{r_2})$ is the complementary design of $(V, \mathcal B_{r_1})$, i.e.,
 $\mathcal B_{r_2}=\{V\backslash B\mid B\in \mathcal B_{r_1}\}$.
 \item  The last column for each $n$ denotes the number of non-isomorphic symmetric designs.
 `` $\times$ '' indicates the non-existence,  `` $?$ '' indicates that existence or non-existence is unknown.
 \item "$\sharp$" denotes the number of non-isomorphic designs.
       The information is basically from the Appendix-Tables A, B and C in \cite{Beth-J-L-1985} and Table 1.35 in \cite{Colbourn-D-2007}.
  \end{enumerate}

 \noindent
 {\bf{Case 2: }} $r_1+r_2=n$ and $w_{r_1}\neq w_{r_2}$.\\
 In the table below, we give the possible values of the pair $(\lambda_3^{(r_1)}(i_1,i_2,i_3)$, $\lambda_3^{(r_2)}(i_1,i_2,i_3))$
 for 3-point subset $\{i_1,i_2,i_3\}\subset V$  and ``$\times$'' indicates the non-existence of tight relative $3$-design with the corresponding parameters.

\begin{spacing}{1.1}
\begin{center}
{\footnotesize
\begin{tabular}{ccc}
 $\begin{array}{|c|c|c|c|c|c|}
\hline
n & r_1 & \lambda_2^{(r_1)} & \frac{w_{r_2}}{w_{r_1}} &(\lambda_3^{(r_1)},\lambda_3^{(r_2)}) & \\ \hline
37 & 9 & 2 & \frac{2}{7} & (0,17),(2,10) & \times \\ \hline
& & &  \frac{1}{6} & (0,18),(1,12),(2,6) & \times \\ \hline
& & &  \frac{2}{17} & (0,19),(2,2)  & \times\\ \hline
& & &  \frac{1}{11} & (0,20),(1,9)  & \times\\ \hline \hline
56 & 11 & 2 & \frac{1}{4} & (0,30),(1,26),(2,22)  &   \\ \hline
& & & \frac{1}{7} & (0,31),(1,24),(2,17)  &   \\ \hline
& & & \frac{1}{10} & (0,32),(1,22),(2,12)  &   \\ \hline
& & & \frac{1}{13} & (0,33),(1,20),(2,7)  &   \\ \hline
& & & \frac{1}{16} & (0,34),(1,18),(2,2)  &   \\ \hline
& & & \frac{1}{19} & (0,35),(1,16)  &   \\ \hline
& & & \frac{1}{22} & (0,36),(1,14)  &   \\ \hline
& & & \frac{1}{10} & (0,32),(1,22),(2,12)  &   \\ \hline \hline
66 & 26 & 10 & \frac{1}{13} & (3,24),(4,11)  &   \\ \hline
& & & \frac{1}{9} & (3,21),(4,12),(5,3)  &   \\ \hline
& & & \frac{1}{5} & (i,33-5i),i=2,\ldots,6  &   \\ \hline
& & & \frac{5}{13} & (0,24),(5,11)  &   \\ \hline
& & & \frac{5}{9} & (0,21),(5,12),(10,3)  &   \\ \hline
& & & 5 & (0,15),(5,14),(10,13)  &   \\ \hline
& & & \frac{3}{19} & (3,19),(6,0)  &   \\ \hline
& & & \frac{3}{11} & (3,17),(6,6)  &   \\ \hline
& & & \frac{3}{7} & (3i,23-7i),i=0,1,2,3  &   \\ \hline
& & & \frac{7}{19} & (2,19),(9,0)  &   \\ \hline
& & & \frac{7}{15} & (2,18),(9,3)  &   \\ \hline
& & & \frac{7}{11} & (2,17),(9,6)  &   \\ \hline
& & & \frac{7}{3} & (2,15),(9,12)  &   \\ \hline
\end{array}$
&
$\begin{array}{|c|c|c|c|c|c|}
\hline
n & r_1 & \lambda_2^{(r_1)} & \frac{w_{r_2}}{w_{r_1}} &(\lambda_3^{(r_1)},\lambda_3^{(r_2)}) & \\ \hline
70 & 24 & 8 & \frac{3}{20} & (1,30),(4,10)  &   \\ \hline
& & & \frac{2}{19} & (2,25),(4,6)  &   \\ \hline
& & & \frac{1}{18} & (2,30),(3,12)  &   \\ \hline
& & & \frac{5}{22} & (2,22),(7,0)  &   \\ \hline \hline
71 & 15 & 3 & \frac{2}{25} & (1,29),(3,4)  &   \\ \hline
 &  & & \frac{1}{24} & (1,24),(2,0)  &   \\ \hline
 & 21 & 6 & \frac{3}{26} & (1,30),(4,4)  &   \\ \hline
 &  &  & \frac{5}{28} & (1,28),(6,0)  &   \\ \hline \hline
78 & 22 & 6 & \frac{2}{21} & (2,24),(4,3)  &   \\ \hline
& & & \frac{4}{23} & (2,26),(6,3)  &   \\ \hline
& & & \frac{1}{20} & (1,40),(2,20),(3,0))  &   \\ \hline
& & & \frac{6}{25} & (0,35),(6,10)  &   \\ \hline
& & & \frac{5}{24} & (0,36),(5,12)  &   \\ \hline
& & & \frac{3}{22} & (0,40),(3,18)  &   \\ \hline \hline
79 & 13 & 2 & \frac{2}{9} & (0,47),(2,38)  &   \\ \hline
& & & \frac{1}{8} & (0,48),(1,40),(2,32)  &   \\ \hline
& & & \frac{2}{23} & (0,49),(2,26)  &   \\ \hline
& & & \frac{1}{15} & (0,50),(1,35),(2,20)  &   \\ \hline
& & & \frac{2}{37} & (0,51),(2,14)  &   \\ \hline
& & & \frac{1}{22} & (0,52),(1,30),(2,8)  &   \\ \hline
& & & \frac{2}{51} & (0,53),(2,2)  &   \\ \hline
& & & \frac{1}{29} & (0,54),(1,25)  &   \\ \hline \hline
96 & 20 & 4 & \frac{4}{51} & (0,57),(4,6)  &   \\ \hline
 &  &  & \frac{3}{50} & (0,60),(3,10) &   \\ \hline
\end{array}$\\
\\
\end{tabular}}
	{Table 2}
\end{center}
\end{spacing}
\vspace{-0.2in}
 \begin{pro}
 There is no tight relative 3-design on two shells with $n=37, r_1=9, r_2=28$ and $w_{r_1}\neq w_{r_2}$.
  \end{pro}
 {\bf Proof. ~}
It is known that there exist exactly four 2-$(37,9,2)$ symmetric designs.
For all the four $2$-$(37,9,2)$ designs we proved by computation that the $\lambda_3$-sequences equal to $(4662*0, 3108*1)$.
 Hence $\lambda_3$-sequences of four 2-$(37,28,21)$ symmetric designs equal to $(3108*15, 4662*16)$.
This implies that $\lambda_3^{(r_2)}(i_1,i_2,i_3)=15$ or 16 for any 3 points $i_1,i_2,i_3 \in V$.
Hence it is impossible to have tight relative 3-design with this parameter.
\hfill\qed

 \begin{pro}
 The tight relative 3-designs on two shells with $n=66, r_1=26, r_2=40$ and $w_{r_1}\neq w_{r_2}$ corresponding to the $14$  known symmetric designs in the home page of Ted Spence $($see \cite{Spence}$)$ do not exist.
  \end{pro}
  {\bf Proof. ~}
So far, the number of the known non-isomorphic  $2$-$(66,26,10)$ designs is 14.
They give 14 different types of  $\lambda_3$-sequence.
$${\footnotesize
\begin{array}{llllllllllll}
(110*0, & 825*1, & 3850*2, & 9900*3, & 23100*4, & 6875*5, & 1100*6), \\
(110*0, & 475*1, & 3950*2, & 11200*3, & 21950*4, & 6725*5, & 1300*6, & 50*8), \\
(60*0, & 425*1, & 3650*2, & 12775*3, & 20125*4, & 7525*5, & 1025*6, & 150*7, & 25*9), \\
(60*0, & 450*1, & 3675*2, & 12775*3, & 19775*4, & 7950*5, & 950*6, & 100*7, & 25*9), \\
(60*0, & 575*1, & 3875*2, & 11850*3, & 20175*4, & 8475*5, & 625*6, & 100*7, & 25*8), \\
(85*0, & 400*1, & 3650*2, & 12650*3, & 20400*4, & 7300*5, & 1125*6, & 125*7, & 25*9), \\
(85*0, & 425*1, & 3600*2, & 12775*3, & 20050*4, & 7725*5, & 925*6, & 150*7, & 25*9), \\
(85*0, & 600*1, & 3725*2, & 11900*3, & 20300*4, & 8400*5, & 625*6, & 100*7, & 25*8), \\
(110*0, & 325*1, & 3700*2, & 12450*3, & 21000*4, & 6775*5, & 1250*6, & 50*7, & 100*8), \\
(110*0, & 525*1, & 3625*2, & 11675*3, & 21825*4, & 6725*5, & 1125*6, & 75*7, & 75*8), \\
(110*0, & 575*1, & 3775*2, & 11175*3, & 22050*4, & 6825*5, & 1175*6, & 25*7, & 50*8), \\
(110*0, & 600*1, & 3750*2, & 11050*3, & 22300*4, & 6700*5, & 1150*6, & 50*7, & 50*8), \\
(135*0, & 475*1, & 3775*2, & 11425*3, & 21925*4, & 6675*5, & 1275*6, & 25*7, & 50*8), \\
(135*0, & 675*1, & 3700*2, & 10650*3, & 22750*4, & 6625*5, & 1150*6, & 50*7, & 25*8). \\
\end{array}}
$$
However, the list $(\lambda_3^{(r_1)},\lambda_3^{(r_2)})$ for $n=66$ is given in Table 2.
Hence it is impossible to have tight relative 3-designs for $n=66$ with non-constant weight if $(V,\mathcal B_{r_1})$ is the one of the fourteen $2$-$(66,26,10)$
designs.\\

\noindent
{\bf{Case 3: }} $r_1+r_2\neq n$ and $w_{r_1}= w_{r_2}$.\\
The following is the table of the feasible parameters for $n\leq 100$.

\begin{center}
\begin{tabular}{ccc}
$
\begin{array}{|c|c|c||c|c||c|c|}
 \hline
n & r_1 & \lambda_2^{(r_1)} &r_2 & \lambda_2^{(r_2)} &
(\lambda_3^{(r_1)},\lambda_3^{(r_2)})&\\ \hline \hline
31 & 6 & 1 & 16 & 8 &(i,4-i),i=0,1&\times \\ \hline
31 & 15 & 7 & 25 & 20 & (i,19-i),i=0,\ldots,7&\times \\ \hline \hline
85 & 21 & 5 & 49 & 28 & (i,17-i), i=0,\ldots,5& \\
\hline85 & 36 & 15 & 64 & 48 & (i,42-i), i=0,\ldots,15& \\\hline
\end{array}$\\
\\
	{Table 3}
\end{tabular}
\end{center}

\begin{pro}
There is no tight relative $3$-design with $n=31$, $r_1+r_2\neq 31$ and $w_{r_1}=w_{r_2}$.
\label{pro:5.3}
\end{pro}
{\bf Proof. ~}
(1) Non-existence for $n=31$, $r_1=6$, $r_2=16$.\\
\eqref{equ:5.2} implies $\lambda^{(r_1)}_3+\lambda^{(r_2)}_3=4$.
The $\lambda_3$-sequence of the 2-$(31,6,1)$ design is $(3875*0, 620*1)$.
Then we must have $\lambda^{\mathcal B_{r_2}}_3(i_1,i_2,i_3)=4$ or $3$ for any 3-point subset $\{i_1,i_2,i_3\}\subset V$.
 The complementary design $(V,\mathcal B_{r_2}^c)$ is a $2$-$(31,15, 7)$ design and \eqref{equ:4.11} implies
 $\lambda^{\mathcal B_{r_2}^c}_3(i_1,i_2,i_3) +\lambda^{\mathcal B_{r_2}}_3(i_1,i_2,i_3)=7$.
 Hence $\lambda^{\mathcal B_{r_2}^c}_3(i_1,i_2,i_3)=3$ or $4$.
 Hence the $\lambda_3$-sequence of $(V,\mathcal B_{r_2}^c)$ must be of the form $(a_3\ast 3, a_4\ast 4)$.
On the other hand there are 10374196953 non-isomorphic 2-$(31,15,7)$ designs, and Brendan McKay got all the $\lambda_3$-sequences with length at most 4 as listed below.
We can easily see that there is no $\lambda_3$-sequence satisfying this condition.
Therefore there is no tight relative 3-design with this parameter.

 {\footnotesize
 	\begin{eqnarray}
 	\begin{array}{|l|l|}
 	No.& \\ \hline
 	1 & (4340*3, \  155*7) \  \\
 	2 & (930*2, \  2015*3, \ 1550*4) \ \\
 	3 & (64*0, \  3892*3, \ 448*4, \ 91*7)  \\
 	4 & (112*0, \ 3556*3, \ 784*4, \ 43*7)  \\
 	5 & (48*1, \ 4196*3, \ 144*5, \ 107*7)  \\
 	6 & (64*1, \ 4148*3, \ 192*5, \ 91*7)  \\
 	7 & (72*1, \ 4124*3, \ 216*5, \ 83*7)  \\
 	8 & (80*1, \ 4100*3, \ 240*5, \ 75*7)  \\
 	9 & (84*1, \ 4088*3, \ 252*5, \ 71*7)  \\
 	10 & (88*1, \ 4076*3, \ 264*5, \ 67*7)  \\
 	11 & (92*1, \ 4064*3, \ 276*5, \ 63*7)  \\
 	12 & (96*1, \ 4052*3, \ 288*5, \ 59*7)  \\
 	13 & (100*1, \ 4040*3, \ 300*5, \ 55*7)  \\
 	14 & (104*1, \ 4028*3, \ 312*5, \ 51*7)  \\ \hline
 	\end{array}
 	&
 	\begin{array}{|l|l|}
 	No.& \\ \hline
 	15 & (108*1, \ 4016*3, \ 324*5, \ 47*7)  \\
 	16 & (112*1, \ 4004*3, \ 336*5, \ 43*7)  \\
 	17 & (116*1, \ 3992*3, \ 348*5, \ 39*7)  \\
 	18 & (120*1, \ 3980*3, \ 360*5, \ 35*7)  \\
 	19 & (124*1, \ 3968*3, \ 372*5, \ 31*7)  \\
 	20 & (128*1, \ 3956*3, \ 384*5, \ 27*7)  \\
 	21 & (132*1, \ 3944*3, \ 396*5, \ 23*7)  \\
 	22 & (136*1, \ 3932*3, \ 408*5, \ 19*7)  \\
 	23 & (140*1, \ 3920*3, \ 420*5, \ 15*7)  \\
 	24 & (144*1, \ 3908*3, \ 432*5, \ 11*7)  \\
 	25 & (148*1, \ 3896*3, \ 444*5, \ 7*7)  \\
 	26 & (840*2, \ 2285*3, \ 1280*4, \ 90*5)  \\
 	27 & (855*2, \ 2240*3, \ 1325*4, \ 75*5) \\ \hline
 	\end{array}\label{equ:5.3}
 	\end{eqnarray}
}

\noindent
(2) Non-existence for $n=31$, $r_1=15$, $r_2=25$.\\
By \eqref{equ:5.2}, for $r_1=15$ and $r_2=25$, we have
$$
\lambda^{(r_1)}_3(i_1,i_2,i_3)+\lambda^{(r_2)}_3(i_1,i_2,i_3)=19.
$$
Since $r_2=25$, $(V,\mathcal B_{r_2}^c)$ is a symmetric $2$-$(31,6,1)$ design.
Hence \eqref{equ:4.11} implies
$$
\lambda^{\mathcal B_{r_2}}_3(i_1,i_2,i_3)=16-\lambda^{\mathcal B_{r_2}^c}_3(i_1,i_2,i_3)=15,~ \mbox{or}~ 16.
$$
Therefore for $r_1=15$, we have
$$
\lambda^{\mathcal B_{r_1}}_3(i_1,i_2,i_3)=3~ \mbox{or}~ 4.
$$
Therefore \eqref{equ:5.3} implies the non-existence of tight relative 3-design of this parameter.
\hfill\qed\\

\noindent
{\bf{Case 4: }}$r_1+r_2\neq n$ and $w_{r_1}\neq w_{r_2}$.\\
\begin{spacing}{1.1}
	$${\footnotesize
	\begin{array}{c||c|c||c|c||c|cc}
	n & r_1 & \lambda_2^{(r_1)} &r_2 & \lambda_2^{(r_2)} & \frac{w_{r_2}}{w_{r_1}} & (\lambda_3^{(r_1)},\lambda_3^{(r_2)})& \\ \hline
	31 & 6 & 1 & 21 & 14  & \frac{1}{6} & (0,10),(1,4) & \times\\ \hline
	& 10 & 3 & 16 & 8  & \frac{1}{5} & (0,8),(1,3) & \times\\ \hline
	& 10 & 3 & 25 & 20  & \frac{1}{5} &(i,20-5i),i=0,1,2,3 & \times\\ \hline
	& 15 & 7 & 21 & 14  & \frac{1}{6} & (3,10),(4,4) & \times\\ \hline
	& 16 & 8 & 21 & 14  & \frac{4}{5} & (0,14),(4,9),(8,4) & \times\\ \hline
	& 21 & 14 & 25 & 20  & \frac{5}{4} & (4,20),(9,16),(14,12) & \times\\ \hline
	& 21 & 14 & 25 & 20  & \frac{4}{9} & (10,14),(14,5) & \times\\ \hline  \hline
	\end{array}}$$
\end{spacing}
	
\begin{spacing}{1.1}
	$${\footnotesize
	\begin{array}{c||c|c||c|c||c|cc}
	n & r_1 & \lambda_2^{(r_1)} &r_2 & \lambda_2^{(r_2)} & \frac{w_{r_2}}{w_{r_1}} & (\lambda_3^{(r_1)},\lambda_3^{(r_2)}) \\ \hline
	61 & 16 & 4 & 36 & 21  & \frac{1}{2} & (i,14-2i),i=1,2,3,4 & \\ \hline
	& 25 & 10 & 45 & 33  & 2 & (2i,26-i),i=0,\ldots,5 & \\ \hline
	&   &   &   &    & \frac{7}{33} & (2,33),(9,0) & \\ \hline
	& 36 & 21 & 45 & 33  & \frac{15}{22} & (6,33),(21,11) & \\ \hline
	&  &  &  &   & \frac{11}{24} & (8,33),(19,9) & \\ \hline
	&  &  &  &  & \frac{5}{27} & (11,30),(16,3) & \\ \hline \hline
	71 & 15 & 3 & 36 & 18  & \frac{1}{2} & (i,10-2i),i=0,1,2,3  & \\ \hline
	& 15 & 3 & 50 & 35  & \frac{1}{10} &(0,30), (1,20),(2,10) & \\ \hline
	& 21 & 6 & 35 & 17  & \frac{3}{4} & (1,9),(4,5) & \\ \hline
	&   &   &   &    & \frac{1}{9} & (1,14),(2,5) & \\ \hline
	& 21 & 6 & 36 & 18  & \frac{2}{5} & (0,13),(2,8),(4,3) & \\ \hline
	&   &   &   &    & 5 & (1,9),(6,8) & \\ \hline
	&   &   &   &    & \frac{1}{14} & (1,18),(2,4) & \\ \hline
	& 21 & 6 & 56 & 44  & \frac{5}{12} & (1,36),(6,24) & \\ \hline
	&   &   &   &    & \frac{3}{44} & (1,44),(4,0) & \\ \hline
	&   &   &   &    & \frac{4}{5} & (2,34),(6,29) & \\ \hline
	&   &   &   &    & \frac{1}{7} & (i,46-7i),i=1,\ldots,6  & \\ \hline
	& 35 & 17 & 50 & 35  & \frac{11}{9} & (0,31),(11,22) & \\ \hline
	& & & & & \frac{17}{16} & (0,32),(17,16) & \\ \hline
	& & & & & \frac{5}{2} & (4,26),(9,24),(14,22) & \\ \hline
	& & & & & \frac{6}{7} & (5,28),(11,21),(17,14) & \\ \hline
	& & & & & \frac{8}{17} & (5,31),(13,14) & \\ \hline
	& & & & & \frac{9}{22} & (5,32),(14,10) & \\ \hline
	& & & & & \frac{7}{12} & (6,28),(13,16) & \\ \hline
	& & & & & \frac{1}{28} & (8,28),(9,0) & \\ \hline
	& & & & & \frac{1}{5} & (i,65-5i),i=6,\ldots,12 & \\ \hline
	& 35 & 17 & 56 & 44  & 2 & (2i-1,39-i),i=1,\ldots,9 & \\ \hline
	& & & & & \frac{17}{20} & (0,44),(17,24) & \\ \hline
	& & & & & \frac{9}{16} & (5,40),(14,24) & \\ \hline
	& & & & & \frac{12}{29} & (5,42),(17,13) & \\ \hline
	& & & & & \frac{7}{15} & (6,39),(13,24) & \\ \hline
	& & & & & \frac{3}{13} & (8,35),(11,22),(14,9) & \\ \hline
	& & & & & \frac{1}{12} & (8,36),(9,24),(10,12) & \\ \hline
	& & & & & \frac{5}{14} & (9,32),(14,18) & \\ \hline \hline	71 & 35 & 17 & 56 & 44  & \frac{1}{35} & (8,39),(9,4) & \\ \hline
	& & & & & \frac{4}{25} & (9,29),(13,4) & \\ \hline
	& 36 & 18 & 50 & 35  & \frac{1}{18} & (9,22),(10,4) & \\ \hline
	& & & & & \frac{17}{7} & (0,28),(17,21) & \\ \hline
	& & & & & \frac{4}{3} & (4i,31-3i),i=0,\ldots,4  & \\ \hline
	& & & & & 9 & (3,25),(12,24) & \\ \hline
	& & & & & \frac{15}{17} & (3,31),(18,14) & \\ \hline
		& & & & & \frac{7}{11} & (4,32),(11,21),(18,10) & \\ \hline
		& & & & & \frac{11}{14} & (6,28),(17,14) & \\ \hline
		& & & & & \frac{3}{8} & (3i,48-8i),i=2,3,4,5  & \\ \hline
	\end{array}}$$
\end{spacing}

\begin{spacing}{1.1}
	$${\footnotesize
	\begin{array}{c||c|c||c|c||c|cc}
	n & r_1 & \lambda_2^{(r_1)} &r_2 & \lambda_2^{(r_2)} & \frac{w_{r_2}}{w_{r_1}} & (\lambda_3^{(r_1)},\lambda_3^{(r_2)}) \\ \hline
	71 & 36 & 18 & 50 & 35  & \frac{10}{19} & (8,26),(18,7) & \\ \hline
		& & & & & \frac{5}{21} & (8,28),(13,7) & \\ \hline
		& & & & & \frac{2}{13} & (8,30),(10,17),(12,4) & \\ \hline
	& 36 & 18 & 56 & 44  & \frac{3}{10} & (3i,64-10i),i=2,\ldots,6 & \\ \hline
	& & & & & \frac{5}{9} & (8,36),(13,27),(18,18) & \\ \hline
	& & & & & \frac{17}{3} & (0,36),(17,33) & \\ \hline
	& & & & & \frac{15}{4} & (3,36),(18,32) & \\ \hline
	& & & & & \frac{9}{7} & (3,39),(12,32) & \\ \hline
	& & & & & \frac{7}{8} & (4,40),(11,32),(18,24) & \\ \hline
	& & & & & \frac{11}{6} & (6,36),(17,30) & \\ \hline
	& & & & & \frac{11}{29} & (6,42),(17,13) & \\ \hline
	& & & & & \frac{4}{21} & (8,39),(12,18) & \\ \hline
	& & & & & \frac{5}{32} & (8,40),(13,8) & \\ \hline
	& & & & & \frac{1}{11} & (i,132-11i),i=8,9,10,11  & \\ \hline
	& 50 & 35 & 56 & 44  & \frac{16}{3} & (0,39),(16,36),(32,33) & \\ \hline
	& & & & & \frac{35}{8} & (0,40),(35,32) & \\ \hline
	& & & & & \frac{28}{11} & (0,44),(28,33) & \\ \hline
	& & & & & 13 & (4,36),(17,35),(30,34) & \\ \hline
	& & & & & \frac{19}{5} & (7,39),(26,34) & \\ \hline
	& & & & & \frac{22}{7} & (10,39),(32,32) & \\ \hline
	& & & & & \frac{3}{2} &(3i+1,50-2i),i=3,\ldots,11 & \\ \hline
	& & & & & \frac{5}{11} &(5i,88-11i),i=4,5,6,7 & \\ \hline
	& & & & & \frac{11}{15} & (21,39),(32,24) & \\ \hline
	& & & & & \frac{7}{20} & (21,44),(28,24),(35,4) & \\ \hline
	& & & & & \frac{9}{29} & (22,42),(31,13) & \\ \hline
	& & & & & \frac{8}{13} & (24,35),(32,22) & \\ \hline
	& & & & & \frac{2}{9} & (2i,144-9i),i=12,13,14,15 & \\ \hline
	& & & & & \frac{1}{16} & (24,40),(25,24),(26,8) & \\ \hline
	& & & & & \frac{3}{25} & (25,29),(28,4) & \\ \hline  \hline
	79 & 27 & 9 & 39 & 19  & \frac{3}{5} & (0,14),(3,9),(6,4)\\ \hline
	& 27 & 9 & 66 & 55  & \frac{1}{22} & (3,44),(4,22)\\ \hline
	& 27 & 9 & 66 & 55  & \frac{4}{11} & (1,51),(5,40),(9,29)\\ \hline
	& 39 & 19 & 52 & 34  & \frac{4}{13} & (7,29),(11,16),(15,3)\\ \hline
	& & & & & \frac{19}{4} & (0,24),(19,20)\\ \hline
	& & & & & \frac{9}{10} & (2,30),(11,20)\\ \hline
	& & & & & \frac{3}{29} & (8,33),(11,4)\\ \hline
	& 39 & 19 & 66 & 55  & \frac{19}{11} & (0,51),(19,40)\\ \hline
	& & & & &\frac{6}{19} & (8,28),(14,9)\\ \hline
	& & & & &\frac{1}{16} & (10,20),(11,4)\\ \hline
	& 40 & 20 & 52 & 34  & \frac{5}{3} &(5i,28-3i),i=0,\ldots,4\\ \hline
	& & & & &\frac{6}{19} & (8,28),(14,9)\\ \hline
	& & & & &\frac{1}{16} & (10,20),(11,4)\\ \hline
	& 40 & 20 & 66 & 55  & \frac{13}{33} & (7,53),(20,20)\\ \hline
	& & & & &\frac{9}{11} & (8,48),(17,37)\\ \hline
	& & & & &\frac{2}{11} &(2i,100-11i),i=5,\ldots,9\\ \hline
	\end{array}}$$
\end{spacing}

\begin{spacing}{1.1}
	$${\footnotesize
\begin{array}{c||c|c||c|c||c|c}
n & r_1 & \lambda_2^{(r_1)} &r_2 & \lambda_2^{(r_2)} & \frac{w_{r_2}}{w_{r_1}} & (\lambda_3^{(r_1)},\lambda_3^{(r_2)}) \\ \hline
79& 52 & 34 & 66 & 55  & \frac{3}{11} & (3i+1,123-11i),i=7,\ldots,11\\ \hline
 & & & & &\frac{24}{11} & (4,54),(28,43)\\ \hline
 & & & & &\frac{17}{11} & (17,49),(34,38)\\ \hline
 & & & & &\frac{10}{11} & (20,48),(30,37)\\ \hline
 & & & & &\frac{2}{33} & (22,47),(24,14)\\ \hline \hline
 85 & 21 & 5 & 57 & 38  & \frac{3}{17} & (1,26),(4,9)\\ \hline
   & 28 & 9 & 49 & 28  & \frac{2}{15} & (2,22),(4,7),(6,12)\\ \hline
   & 28 & 9 & 36 & 15  & \frac{5}{4} & (3,6),(8,2)\\ \hline
   & 28 & 9 & 64 & 48  & \frac{2}{15} & (4,27),(6,12)\\ \hline
  & & & &  & \frac{9}{26} & (0,44),(9,18)\\ \hline
  & 36 & 15 & 57 & 38  & \frac{15}{2} & (0,26),(15,24)\\ \hline
  & & & &  & \frac{7}{12} & (1,34),(8,22),(15,10)\\ \hline
   & & & & & \frac{3}{17} & (6,26),(9,9)\\ \hline
  & 49 & 28 & 57 & 38  & \frac{23}{8} & (2,30),(25,22)\\ \hline
  & & & & & \frac{4}{5} & (4i,45-5i),i=2,\ldots,7\\ \hline
  & 57 & 38 & 64 & 48  & 22 & (0,37),(22,36)\\ \hline
  & & & & &\frac{5}{4} & (4i,45-5i),i=2,\ldots,7\\ \hline
  & & & & &\frac{8}{23} & (22,45),(30,22)\\ \hline
  & & & & &\frac{11}{42} & (22,48),(33,6)\\ \hline
  & & & & &\frac{3}{19} & (25,37),(28,18)\\ \hline
  & & & & &\frac{1}{34} & (25,42),(26,8)\\ \hline \hline
 91 & 10 & 1 & 46 & 23  & \frac{1}{7} & (0,12),(1,5)\\ \hline
  & 10 & 1 & 55 & 33  & \frac{1}{15} & (0,21),(1,6)\\ \hline
  & 36 & 14 & 45 & 22  & \frac{2}{5} & (2,19),(4,14),(6,9),(8,4)\\ \hline
  & 36 & 14 & 81 & 72  & \frac{1}{6} & (i,96-6i),i=4,\ldots,14\\ \hline
  & 45 & 22 & 55 & 33  & \frac{22}{9} & (0,24),(22,15)\\ \hline
  & & & & &\frac{12}{13} & (2,29),(14,16)\\ \hline
  & & & & &\frac{7}{15} & (10,21),(17,6)\\ \hline
  & & & & &\frac{2}{17} & (10,25),(12,8)\\ \hline
  & 45 & 22 & 81 & 72  & \frac{16}{15} & (2,72),(18,57)\\ \hline
 & & & & &\frac{11}{27} & (11,63),(22,36)\\ \hline
 & & & & &7 & (3,65),(10,64),(17,63)\\ \hline
 & & & & &\frac{9}{14} & (8,68),(17,54)\\ \hline
 & & & & &\frac{13}{40} & (8,72),(21,32)\\ \hline
 & & & & &\frac{2}{13} & (2i,133-13i),i=12,\ldots,15\\ \hline
  & 46 & 23 & 55 & 33  & \frac{5}{2} & (13,19),(18,17),(23,15)\\ \hline
  & 46 & 23 & 55 & 33  & \frac{23}{27} & (0,33),(23,6)\\ \hline
 & & & & &\frac{3}{19} & (11,22),(14,3)\\ \hline
 & & & & &\frac{5}{2} & (5i+3,23-2i),i=0,\ldots,4\\ \hline
  & 46 & 23 & 81 & 72  & \frac{1}{38} & (12,40),(13,2)\\ \hline
  & & & & &\frac{8}{37} & (12,61),(20,24)\\ \hline
  & & & & &\frac{5}{12} & (5i+3,84-12i),i=1,2,3,4\\ \hline
  & & & & &\frac{3}{25} & (11,67),(14,42),(17,17)\\ \hline
  & & & & &\frac{12}{11} & (8,67),(20,56)\\ \hline
  & & & & &\frac{11}{62} & (11,66),(22,4)\\ \hline
    \end{array}}$$
\end{spacing}

\begin{spacing}{1.1}
	$${\footnotesize
\begin{array}{c||c|c||c|c||c|c}
n & r_1 & \lambda_2^{(r_1)} &r_2 & \lambda_2^{(r_2)} & \frac{w_{r_2}}{w_{r_1}} & (\lambda_3^{(r_1)},\lambda_3^{(r_2)}) \\ \hline
91 & 55 & 33 & 81 & 72  & \frac{8}{41} & (21,57),(29,16)\\ \hline
 & & & & &\frac{10}{29} & (21,60),(31,31)\\ \hline
 & & & & &\frac{12}{17} & (21,62),(33,45)\\ \hline
 & & & & &\frac{11}{23} & (22,59),(33,36)\\ \hline
 & & & & &\frac{9}{35} & (24,47),(33,12)\\ \hline
 & & & & &\frac{27}{16} & (6,72),(33,56)\\ \hline
 & & & & &\frac{14}{5} & (11,67),(25,62)\\ \hline
 & & & & &\frac{13}{11} & (16,67),(29,56)\\ \hline
 \end{array}}$$
 \begin{center}{Table 4}\end{center}
\end{spacing}
\vspace{-0.1in}
 \begin{pro}
 There is no tight relative $3$-design with the parameters $n=31$, $r_1+r_2\neq n$, $w_{r_1}\neq w_{r_2}$ in the list given above.
 \end{pro}
 {\bf Proof. ~}
There are 151 non-isomorphic 2-$(31,10,3)$ designs and the $\lambda_3$-sequence is one of the following cases.
$${\footnotesize
\begin{array}{llll}
( 1240*0, & 2790*1, & 465*2),& \\
( 1239*0, & 2793*1, & 462*2, & 1*3), \\
( 1238*0, & 2796*1, & 459*2, & 2*3), \\
( 1237*0, & 2799*1, & 456*2, & 3*3), \\
( 1236*0, & 2802*1, & 453*2, & 4*3), \\
( 1235*0, & 2805*1, & 450*2, & 5*3), \\
( 1234*0, & 2808*1, & 447*2, & 6*3), \\
( 1233*0, & 2811*1, & 444*2, & 7*3), \\
( 1232*0, & 2814*1, & 441*2, & 8*3), \\
( 1231*0, & 2817*1, & 438*2, & 9*3), \\
(  1225*0, & 2835*1, & 420*2, & 15*3 ).
\end{array} }
$$
The non-existence of such tight relative 3-designs on two shells is obtained (cf. Section \ref{subsec:further_result}). \hfill\qed

\noindent
 \begin{pro}
 	The tight relative 3-designs on two shells with $n=61, r_1=25, r_2=45$ and $w_{r_1}\neq w_{r_2}$ from the 31 known symmetric designs in the home page of Ted Spence  do not exist.
 \end{pro}
 {\bf Proof. ~}
In the home page of Ted Spence, there are 31 non-isomorphic 2-$(61,25,10)$ designs.
And we can calculate the following $\lambda_3$-sequences for these $2$-$(61,25,10)$ designs.

$$
{\footnotesize
\begin{array}{lllllllll}
(105*0, & 125*1, & 2400*2, & 7650*3, & 17625*4, & 6735*5, & 1200*6, & 150*7), \\
(105*0, & 125*1, & 2325*2, & 7950*3, & 17175*4, & 7035*5, & 1125*6, & 150*7), \\
(105*0, & 125*1, & 2375*2, & 7725*3, & 17550*4, & 6760*5, & 1200*6, & 150*7), \\
(105*0, & 125*1, & 2300*2, & 8025*3, & 17100*4, & 7060*5, & 1125*6, & 150*7), \\
(130*0, & 200*1, & 2300*2, & 7450*3, & 17800*4, & 6760*5, & 1300*6, & 50*7), \\
(130*0, & 200*1, & 2275*2, & 7525*3, & 17725*4, & 6785*5, & 1300*6, & 50*7), \\
(130*0, & 200*1, & 2225*2, & 7750*3, & 17350*4, & 7060*5, & 1225*6, & 50*7), \\
(130*0, & 200*1, & 2200*2, & 7825*3, & 17275*4, & 7085*5, & 1225*6, & 50*7), \\
\end{array} }$$
$$
{\footnotesize
	\begin{array}{lllllllll}
(180*0, & 230*1, & 2800*2, & 5715*3, & 18880*4, & 7405*5, & 705*6, & 50*7, & 25*8), \\
(180*0, & 225*1, & 2720*2, & 5875*3, & 18930*4, & 7170*5, & 815*6, & 50*7, & 25*8), \\
(185*0, & 200*1, & 2640*2, & 6360*3, & 18150*4, & 7680*5, & 700*6, & 60*7, & 15*8), \\
(185*0, & 250*1, & 2595*2, & 6250*3, & 18250*4, & 7740*5, & 645*6, & 60*7, & 15*8), \\
(195*0, & 225*1, & 2700*2, & 5905*3, & 18800*4, & 7280*5, & 845*6, & 30*7, & 10*8), \\
(195*0, & 250*1, & 2645*2, & 5905*3, & 18880*4, & 7215*5, & 860*6, & 30*7, & 10*8), \\
(195*0, & 265*1, & 2720*2, & 5695*3, & 18910*4, & 7410*5, & 755*6, & 30*7, & 10*8), \\
(195*0, & 230*1, & 2675*2, & 5955*3, & 18750*4, & 7305*5, & 840*6, & 30*7, & 10*8), \\
(195*0, & 200*1, & 2730*2, & 6005*3, & 18570*4, & 7445*5, & 805*6, & 30*7, & 10*8), \\
(195*0, & 230*1, & 2585*2, & 6135*3, & 18750*4, & 7125*5, & 930*6, & 30*7, & 10*8), \\
(195*0, & 205*1, & 2770*2, & 5855*3, & 18730*4, & 7390*5, & 805*6, & 30*7, & 10*8), \\
(200*0, & 270*1, & 2745*2, & 5625*3, & 18845*4, & 7595*5, & 670*6, & 30*7, & 10*8), \\
(200*0, & 275*1, & 2590*2, & 5895*3, & 18915*4, & 7240*5, & 835*6, & 30*7, & 10*8).
\end{array} }$$
So, we can see the non-existence of such tight relative 3-designs on two shells, if $(V,\mathcal B_{r_1})$ is one of the 31 known 2-$(61,25,10)$ designs listed in the home page of Ted Spence.
\hfill\qed

%%%%
\begin{theo}[{\cite{Beth-J-L-1985}}]\label{theo:5.4}%[Lemma~3.8]
Let $(V,\mathcal{B})$ be a symmetric $2$-$(v,k,\lambda)$ design.
Then $(k-\lambda)^{v-1}$ is a square.
\end{theo}
Theorem~\ref{theo:5.4} implies the non-existence of 2-$(n,r_1,\lambda_2^{(r_1)})$ design with $n$ even in Table 1.
%%%%%
\begin{theo}[{\cite{Beth-J-L-1985}}]\label{theo:5.5}%[Theorem~4.6]
Let $v$ be odd and assume the existence of a symmetric $2$-$(v,k,\lambda)$ design. Then the diophantine equation
$$ x^2 = (k-\lambda)y^2+(-1)^{(v-1)/2}\lambda z^2 $$
has a non-trivial solution in integers.
\end{theo}
%%%%%
Theorem~\ref{theo:5.5} implies the non-existence of 2-$(n,r_1,\lambda_2^{(r_1)})$ design with $n$ odd in Table 1.

%%%%%%%%%%%%%%%%%%%%%%%%%%%%%%%%%%%%%%%%%%%%%%%%

\subsection{Further results}\label{subsec:further_result}
 The discussion in the present paper led to the following Problems and Conjectures. \\

\noindent
{\bf Problem 1.}
If there is any tight relative $3$-design $(Y,w)$ on two shells $X_{r_1}\cup X_{r_2}$ in $H(n,2)$ with constant weight and $r_1+r_2=n$, then is it true that the corresponding symmetric $2$-$(n,r_1,\lambda^{(r_1)}_2)$ design $(V, \mathcal B_{r_1})$ and $2$-$(n,r_2,\lambda^{(r_2)}_2)$ design $(V, \mathcal B_{r_2})$ are complementary designs with each other?
\\

\noindent
{\bf Conjecture 1.}
 Problem 1 is affirmative.\\

\noindent
{\bf Remark.} We note that the same problem is also formulated for tight relative $t=2e+1$ designs in $H(n,2)$, with  $(V, \mathcal B)$  as tight combinatorial $2e$-designs.\\

We can also re-phrase this problem as follows. (This treat the case when $r_1+r_2=n$.) \\

\noindent
{\bf Problem 2}. Are there two symmetric $2$-$(n,k,\lambda)$ designs such that $(V, \mathcal B_1)$ and $(V, \mathcal B_2)$ are different as designs (although they may be or may not be isomorphic as designs) such that their (unordered) $\lambda_3$-sequences coincide?\\

\noindent
{\bf Conjecture 2}. No such two symmetric designs exist satisfying the condition of
Problem 2.\\

\noindent
{\bf Remark.} No such two symmetric designs are known.
We note that the same problem is also formulated for tight relative $(2e+1)$-designs $Y=Y_{r_1} \cup Y_{r_2}$ on two shells of $H(n,2)$, corresponding to tight combinatorial $2e$-designs $(V,\mathcal B_{r_1})$ and $(V,\mathcal B_{r_2})$.
Namely, there are no two tight combinatorial $2e$-designs with the same $\lambda_{2e+1}$-sequence.
 (So far,  no such two symmetric designs are known, and we may conjecture that such examples may not exist.) \\

Here we record some developments on these two problems(and on two conjectures).
Note that if Conjecture 1 holds then Conjecture 2 also holds.

\begin{enumerate}
\item First we proved Conjecture 1 for $n\leq 16$ by adhoc arguments.

Then we proposed these two conjectures in our seminar.
Then our undergraduate students Zongchen Chen and Da Zhao responded, providing the following results.
Their results will be published as an independent paper \cite{Chen-Zhao-2015}.
\begin{enumerate}[(i)]
\item  Conjecture 2 is true for a symmetric $2$-$(n, k,\lambda)$ design if $\lambda=1$, or $2$.
\item  Conjecture 2 is true for a symmetric $2$-$(n, k,\lambda)$ design if $\lambda=3$, provided $k \geq 17$.
\item  Conjecture 1 is true for $2$-$(19,9,4)$ and $2$-$(23,11,5)$ designs.
\end{enumerate}
 To prove (iii) for $2$-$(19,9,4)$ designs, we needed the information on the incidence metrics of all the four $2$-$(19,9,4)$ designs in the home page of Ted Spence.
As for $2$-$(23,11,5)$ designs, we used all the 1106 incidence matrices provided with the curtesy of Ted Spence; note that only 197 of those with non-trivial automorphism group are listed in his home page.
Then by calculating the $\lambda_3$-sequences (and the automorphism group of the 3-subset multiplicity graph defined in the following) of all those symmetric designs and proved that Conjecture 1 is true for $n=19$ and $23$.
(Later, Chen and Zhao also succeeded in proving Conjecture 1 for $n=27$, by obtaining the incidence matrices of all the 208310 of 2-$(27,13,6)$ designs from the list of Hadamard matrices of order 28.)
Here we summarize the main result (techniques) of Chen and Zhao \cite{Chen-Zhao-2015} below.\\

Let $(V,\mathcal B)$ be a $2$-$(n,k,\lambda)$ design and $\lambda_3(i_1,i_2,i_3)=|\{B\in \mathcal B\mid \{i_1,i_2,i_3\}\subset B\}|$.
Let $\Gamma$ be the Johnson graph $J(n,3)$.
The vertex set of $\Gamma$ is the set of all the 3-point subset of $V$ denoted by ${V\choose 3}$.
We assign weight $\lambda_3(i_1,i_2,i_3)$ for each vertex $\{i_1,i_2,i_3\}\in {V\choose 3}$.
We call $(\Gamma,\lambda_3)$ the 3-subset multiplicity graph of the symmetric design $(V,\mathcal B)$.
Let $\mbox{Aut}(\Gamma,\lambda_3)$ be the subgroup of the automorphism group of $\Gamma$ preserving the weight $\lambda_3$.
\begin{enumerate}[(a)]
\item If two distinct $2$-$(n,k,\lambda)$ designs $(V,\mathcal B_1)$ and $(V,\mathcal B_2)$ are isomorphic, then their 3-subset multiplicity graphs are isomorphic.
(So if the corresponding 3-subset multiplicity graphs are not isomorphic then the two designs are not isomorphic.)
\item Let $(V, \mathcal B)$ be a $2$-$(n,k,\lambda)$ design with $n\geq 7$. Let $\mbox{Aut}(V,\mathcal B)$ be the automorphism group of $(V, \mathcal B)$.
Then $\mbox{Aut}(\Gamma,\lambda_3)\geq \mbox{Aut}(V,\mathcal B)$ holds.
Moreover the following (i) and (ii) hold.
\begin{enumerate}[(i)]
\item  If $\mbox{Aut}(\Gamma,\lambda_3)> \mbox{Aut}(V,\mathcal B)$, then there exists another $2$-$(n,k,\lambda)$ design $(V,\widetilde{\mathcal B})$ isomorphic to $(V,\mathcal B)$ having the same 3-subset multiplicity graph.
\item If $\mbox{Aut}(\Gamma,\lambda_3)= \mbox{Aut}(V,\mathcal B)$, then there are no other design with the same3-subset multiplicity graph.
\end{enumerate}
(If $2$-$(n,k,\lambda)$ design $(V,\mathcal B)$ is given explicitly and if $n$ is small, then one can determine all the 3-subset multiplicity graphs and their automorphism groups by computer.)
\end{enumerate}

\item Then we obtained the following general result  (without needing incidence matrices of each design).
Conjecture 2 is true for any $2$-$(4u-1,2u-1,u-1)$ Hadamard design.
(See Theorem \ref{theo:3.6} for more detail.)
\item We proved Conjecture 2 is true for $2$-$(31,10,3)$ and $2$-$(31,15,7)$.
As for $2$-$(31,10,3)$ designs we calculated $\lambda_3$-sequences using all the 151 incidence matrices of $2$-$(31,10,3)$ designs from Spence \cite{Spence-1992} (provided with the courtesy of Ted Spence;
note that only 44 of those with non-trivial automorphism group are listed in his home page \cite{Spence}).
As for $2$-$(31,15,7)$ designs we used $\lambda_3$-sequence with the length at most 4 of all the 10,374,196,953 symmetric $2$-$(31,15,7)$ designs (coming from Hadamard matrices of order 32 \cite{Kharaghani-T-2013}) calculated by Brendan McKay.
 (See the proof of Proposition \ref{pro:5.3} for the detail.)
 \item
 If we could get incidence matrices of all the 78 of $2$-$(25,9,3)$ designs (those 40 of them with non-trivial automorphism group are listed in the home page of Ted Spence), by calculating their $\lambda_3$-sequences, we expect we can show Conjecture 1 is true for $n=25$.
 In the meantime, Chen and Zhao obtained the incidence matrices of all the 2-(25,9,3) designs from Denniston \cite{Denniston-1982}, and then succeeded in proving Conjecture 1 for $n=25$.
 \item Combining all the results mentioned above, we can conclude that Conjecture 1 is true for all $n\leq 35$.
 \end{enumerate}

\subsection{Tight relative $5$-designs}

\begin{theo}
Let $(Y,w)$ be a tight relative $5$-design of $H(n,2)$.
Assume $w$ is constant $w_{r_\nu}$ on each   $Y_{r_\nu}$ for $\nu=1,2$.
Then $n=23$ and $(V,\mathcal B_{r_1})$ and  $(V,\mathcal B_{r_2})$ are combinatorial tight $4$-$(23,7,1)$ design and $4$-$(23,16,52)$ design respectively. Moreover $w_{r_1}=w_{r_2}$ holds and $(V,\mathcal B_{r_1})$ is the complementary design of $(V,\mathcal B_{r_2})$.
\end{theo}
{\bf Proof.}
 By Theorem \ref{theo:3.3}, $(V,\mathcal B_{r_1})$ and $(V,\mathcal B_{r_2})$ are combinatorial tight 4-designs.
 It is well known that there are two combinatorial tight 4-designs, $4$-$(23,7,1)$ design and $4$-$(23,16,52)$ design.
Hence $r_1=7$, $r_2=16$ and $|Y_{r_1}|=|Y_{r_2}|=253$ hold.
Then \eqref{equ:3.1} implies
\begin{eqnarray}
w_{r_1}\lambda^{(r_1)}_5(i_1,i_ 2, i_3,i_ 4, i_5)+w_{r_2}\lambda^{(r_2)}_5(i_1,i_ 2, i_3,i_ 4, i_5)=\frac{3}{19}w_{r_1}+\frac{624}{19}w_{r_2}.
\label{equ:5.4}\end{eqnarray}
Theorem \ref{theo:3.5} implies if $(V,\mathcal B_{r_1})$ and $(V,\mathcal B_{r_2})$ are complementary designs with each other, then the corresponding points set $Y_{r_1}\subset X_{r_1}$ and $Y_{r_2}\subset X_{r_2}$ in $H(n,2)$ gives a tight relative $5$-design$(Y,w)$, with constant weight $w (\equiv 1)$ where $Y=Y_{r_1}\cup Y_{r_2}$. Then in this case \eqref{equ:5.4} implies
$$\lambda^{(r_1)}_5(i_1,i_ 2, i_3,i_ 4, i_5)+\lambda^{(r_2)}_5(i_1,i_ 2, i_3,i_ 4, i_5)=\frac{3}{19}+\frac{624}{19}=33.$$
Since $\lambda^{(r_1)}_4=1$, we have $\lambda^{(r_1)}_5(i_1,i_ 2, i_3,i_ 4, i_5)=0$ or $1$.
 Therefore we have
$\lambda^{(r_2)}_5(i_1,i_ 2, i_3,i_ 4, i_5)=33$ or $32$.
Even if they are not complementary to each other, since there exist 5-point sets $\{a_1,\ldots,a_5\}$ and $\{b_1,\ldots,b_5\}$ in $V$ satisfying
$\lambda^{(r_1)}_5(a_1,\ldots,a_5)=1$ and $\lambda^{(r_1)}_5(b_1,\ldots,b_5)=0$,
\eqref{equ:5.4} implies
$$\lambda^{(r_2)}_5(b_1,b_ 2, b_3,b_ 4, b_5)-\lambda^{(r_2)}_5(a_1,a_ 2, a_3,a_ 4, a_5)=\frac{w_{r_1}}{w_{r_2}}.$$
Therefore we must have $w_{r_2}=w_{r_1}$.
Next assume there exists another $4$-$(23,7, 1)$ design $(V,\mathcal B^\prime_{r_1})$ and $\mathcal B^\prime_{r_1}\cup \mathcal B_{r_2}$ corresponds to the tight relative $5$-design with constant weight.
Let $i_1,i_2,i_3,i_4,i_5$ be any 5 points with $\lambda^{(r_1)}_5(i_1,i_2,i_3,i_4,i_5)=1$.
Then there exists a unique block $B\in \mathcal B_{r_1}$ and $B^\prime\in \mathcal B^\prime_{r_1}$containing $\{i_1,i_2,i_3,i_4,i_5\}$.
Then $\{i_1,i_2,i_3,i_4\}\subset B^\prime,B$.
Let $z^\prime\in B^\prime\backslash \{i_1,i_2,i_3,i_4\}$. Then$\{i_1,i_2,i_3,i_4,z^\prime\}\subset B^\prime$.
Therefore $\lambda^{(r_1)}_5(i_1,i_2,i_3,i_4,z^\prime)=1$.
This implies there exists a unique block in $\mathcal B_{r_1}$ containing $\{i_1,i_2,i_3,i_4,z^\prime\}$.
Since $B$ is the unique block containing $\{i_1,i_2,i_3,i_4\}$, $B$ must also contain $z^\prime$.
Hence we must have $B^\prime\subset B$.
 Since $|B'|=|B|=r_1$, we must have $B'=B$.
 Thus $\mathcal B'_{r_1}$ and $\mathcal B_{r_1}$ coincide.\hfill\qed
\\

\noindent
{\bf Remark ~}
For tight relative $(2e+1)$-designs on two shells in $H(n,2)$ with $e\geq 3$, in view of the non-existence results of $($combinatorial$)$ tight $2e$-designs in $\cite{Enomoto-I-N-1979}, \cite{Peterson-1977}, \cite{Bannai-1977}, \cite{Dukes-S-2013}, \cite{Xiang-2015}$,
we can see that there are no tight relative $(2e+1)$-designs on two shells in $H(n,2)$, if $3\leq e\leq 9$.
And moreover, there are only finitely many tight relative $(2e+1)$-designs on two shells for any fixed $e\geq 10.$

\section{Tight relative 4-designs}\label{sec:4-design}
In this section, we consider tight relative $4$-design on two shells.
For any distinct four points $i_1,\ldots,i_4\in V$, formulas \eqref{equ:4.6} and  \eqref{equ:4.7} given in Proposition \ref{pro:4.2} for $t$ = 4 give the following results.
\begin{eqnarray}
&&p_{\mathcal B_{Y_{r_\nu}}}(\ell;i_1,\ldots,i_s)=\frac{{n-s\choose r_\nu-\ell}}{{n\choose r_\nu}}N_{r_\nu}\qquad\mbox{for $s=1\ldots,3$ and $0\leq \ell\leq s$},\nonumber\\
&&p_{\mathcal B_{Y_{r_\nu}}}(\ell;i_1,\ldots,i_4)=\frac{N_{r_\nu}}{{n\choose r_\nu}}\left\{{n-4\choose r_\nu-\ell}-(-1)^{\ell}{n-4\choose r_\nu-4}\right\}+(-1)^{\ell}\lambda^{(r_\nu)}_4(i_1,\ldots,i_4). \nonumber
\end{eqnarray}
for $0\leq\ell\leq 3$.\\
For $t=4$, the equation \eqref{equ:3.1} is equivalent to the following condition:
$$
\sum_{\nu=1}^2w_{r_\nu}\lambda^{(r_\nu)}_4(i_1,\ldots,i_4)=\sum_{\nu=1}^2 N_{r_\nu}w_{r_\nu}\prod_{j=0}^4\frac{r_\nu-j}{n-j}
$$
for any distinct four points $i_1,\ldots,i_4\in V$.
Also by definition of tight relative 4-design, we must have
$|Y|=\dim(L_2(X_{r_1}\cup X_{r_2}))+\dim(L_1(X_{r_1}\cup X_{r_2}))=\mbox{rank}(E_2)+\mbox{rank}(E_1)=\frac{n(n+1)}{2}$.
We search for the parameters $n,r_1,r_2,N_{r_1},N_{r_2}$ which satisfy all the integral conditions.
Then for each feasible parameter, we investigate whether such combinatorial 3-$(n,r_\nu,\lambda_3^{(r_\nu)})$ exists or not.
The following is the list of feasible parameters of tight relative 4-designs for $n \leq 50$.

$${\small
\begin{array}{ccc|cc|cc|c|c}
n & r_1 & r_2 &N_{r_1} & N_{r_2} & \lambda_3^{(r_1)} & \lambda_3^{(r_2)} &   & \\ \hline
11 & 5 & 6 & 33 & 33 & 2 & 4 &  3\text{-}(11,5,2) & \times \\
16 & 6 & 7 & 56 & 80 & 2 & 5 & 3\text{-}(16,6,2)^{[a]} & \times\\
16 & 6 & 9 & 56 & 80 & 2 & 12 & 3\text{-}(16,6,2) & \times\\
16 & 7 & 10 & 80 & 56 & 5 & 12 & 3\text{-}(16,10,12) & \times,  [a]\\
16 & 9 & 10 & 80 & 56 & 12 & 12 & 3\text{-}(16,10,12) & \times,  [a]\\  \hline
22 & 6 & 7 & 77 & 176 & 1 & 4  &  &\circ \\
22 & 6 & 15 & 77 & 176 & 1 & 52 & &\circ \\
22 & 7 & 8 & 88 & 165 & 2 & 6 & 3\text{-}(22,7,2)^{[b]} & \times\\  % \cite{Colbourn-D-2007}
22 & 7 & 14 & 88 & 165 & 2 & 39 & 3\text{-}(22,7,2) & \times\\ % \cite{Colbourn-D-2007}
22 & 7 & 10 & 176 & 77 & 4 & 6& 3\text{-}(22,10,6)^{[c]} & \times\\
22 & 7 & 12 & 176 & 77 & 4 & 11& 3\text{-}(22,12,11) & \times, [c]\\
22 & 7 & 16 & 176 & 77 & 4 & 28&  &\circ\\
22 & 8 & 15 & 165 & 88 & 6 & 26& 3\text{-}(22,15,26) & \times, [b]\\
22 & 10 & 15 & 77 & 176 & 6 & 52& 3\text{-}(22,10,6) & \times\\
22 & 12 & 15 & 77 & 176 & 11 & 52& 3\text{-}(22,12,11) & \times, [c]\\
22 & 14 & 15 & 165 & 88 & 39 & 26& 3\text{-}(22,15,26) & \times, [b]\\
22 & 15 & 16 & 176 & 77 & 52 & 28&  &\circ\\ \hline
37 & 9 & 10 & 185 & 518 & 2 & 8& 3\text{-}(37,9,2)^{[d]} & \times\\
37 & 9 & 27 & 185 & 518 & 2 & 195& 3\text{-}(37,9,2) & \times\\
37 & 9 & 16 & 370 & 333 & 4 & 24 & &\\
37 & 9 & 21 & 370 & 333 & 4 & 57 & &\\
37 & 10 & 28 & 518 & 185 & 8 & 783 & 3\text{-}(37,28,78) & \times,[d]\\
37 & 16 & 28 & 333 & 370 & 24 & 156 & &\\
37 & 21 & 28 & 333 & 370 & 57 & 156 & &\\
37 & 27 & 28 & 518 & 185 & 195 & 78& 3\text{-}(37,28,78) & \times,[d]\\ \hline
41 & 15 & 16 & 328 & 533 & 14 & 28 & &\\
41 & 15 & 25 & 328 & 533 & 14 & 115 & &\\
41 & 16 & 26 & 533 & 328 & 28 & 80 & &\\
41 & 25 & 26 & 533 & 328 & 115 & 80 & &\\ \hline
46 & 10 & 11 & 253 & 828 & 2 & 9& 3\text{-}(46,10,2)^{[f]} & \times\\
46 & 10 & 35 & 253 & 828 & 2 & 357 &3\text{-}(46,10,2) & \times\\
46 & 11 & 36 & 828 & 253 & 9 & 119 &3\text{-}(46,36,119) & \times,[f]\\
46 & 35 & 36 & 828 & 253 & 357 & 119  &3\text{-}(46,36,119) & \times, [f]\\ \hline
\end{array}
}$$
{\bf Remark}
\begin{enumerate}
	\item The last column denotes the existence and non-existence of 3-$(n,r_{\nu},\lambda_3^{(r_{\nu})})$ design.
\item In the last column, the notation ``$[a]$" denotes that the corresponding $3$-$(n,k,\lambda)$ design is the complementary design of  $3$-$(v,k,\lambda)^{[a]}$ design.
\end{enumerate}

Driessen \cite{Driessen-1978} gave the following condition for the existence of some special $3$-$(n,k,2)$ designs.
\begin{theo}[{\cite{Driessen-1978}}]\label{theo:6.1}%[Theorem 8.4]
	A $3$-$(\binom{u}{2}+u+1,u+1,2)$ design can only exist in one of the following cases:
	\begin{enumerate}
\item $u \equiv 2 \mod 48$ and for every odd prime $p$ and $\alpha$ with $p^{\alpha} \parallel u$ one has $\alpha$ is even or $p \equiv 1,3,9,11 \mod 16$.
\item $u \equiv 14 \mod 48$ and for every odd prime $p$ and $\alpha$ with $p^{\alpha} \parallel u$ one has $\alpha$ is even or $p \equiv 1,7,9,15 \mod 16$.
	\end{enumerate}
\end{theo}
\noindent
Note that the complementary design of a $t$-$(n,k,\lambda)$ design is $t$-$(n,n-k,\mu)$ design with $\mu = \lambda \binom{n-t}{k}/\binom{n-t}{k-t}$.
We consider the complementary design of 3-$(\binom{u+1}{2}+1,u+1,2)$ design, i.e., 3-design with the following parameters.
\[\Big(\binom{u+1}{2}+1,\binom{u}{2},\frac{1}{4}(u^2-u-4)(u-2)\Big).\]
Theorem \ref{theo:6.1} implies the non-existence of some $3$-$(n,k,\lambda)$ designs in the following table.
\[
\begin{array}{c|c|c}
(V,\mathcal B)& (V,\mathcal B^c) & N\\ \hline
3\text{-}(11,5,2) & 3\text{-}(11,6,4) & 33\\
3\text{-}(16,6,2) &3\text{-}(16,10,12) &56\\
3\text{-}(22,7,2) &3\text{-}(22,15,26) &88\\
3\text{-}(37,9,2) &3\text{-}(37,28,78) &185\\
3\text{-}(46,10,2) & 3\text{-}(46,36,119) & 253\\ \hline
\end{array}
\]
Using \eqref{equ:3.1}, we obtain the following lemma.
\begin{lem}\label{lem:extension}
	Let $(V,\mathcal B_r)$ be a $(t-1)$-$(n,r,\lambda)$ design and $(V,\mathcal B_{r+1})$ a $(t-1)$-$(n,r+1,\lambda^{\prime})$ design with $\lambda^{\prime}=\lambda\frac{n-r}{r-t+2}$.
	Let $\infty$ be a point not in $V$ and define $V^+:=V \cup \{\infty\}$ and $\mathcal B_r^+:=\{B \cup \{\infty\} \mid B \in \mathcal B_r\}$.
	If $(V,\mathcal B_r \cup \mathcal B_{r+1})$ is a relative $t$-design with constant weight, then $(V^+,\mathcal B_r^+ \cup \mathcal B_{r+1})$ is a $t$-$(n+1,r+1,\lambda)$ design.
\end{lem}

\begin{theo}\label{pro:6.3}
There exist exactly four tight relative $4$-designs with constant weight when $n=22$ and $(r_1,r_2)=(6,7)$, $(6,15)$, $(7,16)$, $(15,16)$.
%$(r_1,r_2,N_{r_1},N_{r_2})=(6,7,77,176)$, $(6,15,77,176)$, $(7,16,176,77)$, $(15,16,176,77)$.
\end{theo}
{\bf Proof.}
It is proved that
$$
\sum_{\nu=1}^2w_{r_\nu}\lambda^{(r_\nu)}_t(i_1,\ldots,i_t)=\sum_{\nu=1}^2N_{r_\nu}w_{r_\nu}\prod_{j=0}^{t-1}\frac{r_\nu-j}{n-j}.
$$
Putting $w_{r_1}=w_{r_2}$ and $t=4$, for any four distinct points $i_1,i_2,i_3,i_4 \in V$, we obtain that
\begin{equation}
\label{equ:6.1}
\lambda_4^{(r_1)}(i_1,\ldots,i_4)+\lambda_4^{(r_2)}(i_1,\ldots,i_4)=1,33,20,52,
\end{equation}
corresponding to the four cases.
It is known that there is a unique 5-$(24,8,1)$ design which is called the Witt design, so does its derived design $4$-$(23,7,1)$ design $(V^+,\mathcal B)$.
Then we obtain $3$-$(22,6,1)$ design $(V,\mathcal B_{r_1})$ and $3$-$(22,7,4)$ design $(V,\mathcal B_{r_2})$ as the derived design and residual design of $(V^+,\mathcal B)$, as well as their complement designs, i.e., $3$-$(22,16,28)$ design $(V,\mathcal B_{r_1}^c)$ and $3$-$(22,15,52)$ design $(V,\mathcal B_{r_2}^c)$.\\
Using the incidence matrix of Witt design, we can obtain the incidence matrix of $(V,\mathcal B_{r_1})$ and $(V,\mathcal B_{r_2})$. It is not difficult to check \eqref{equ:6.1} is satisfied.
Hence we find four tight relative 4-designs, i.e., $Y_{r_1} \cup Y_{r_2}$, $Y_{r_1}^c \cup Y_{r_2}^c$, $Y_{r_1} \cup Y_{r_2}^c$ and $Y_{r_1}^c\cup Y_{r_2}$.
Here $Y_{r_\nu}$ and $Y_{r_\nu}^c$ correspond to the block sets $\mathcal B_{r_\nu}$ and $\mathcal B_{r_\nu}^c$ ($\nu=1,2$), respectively.\\
We shall prove the uniqueness of each tight relative 4-design.
(Note $4$-$(23,7,1)$ design and 3-$(22,6,1)$ design uniquely exist \cite{Beth-J-L-1985}.)
Let $(V,\mathcal B_6)$ and $(V,\mathcal B_7)$ be derived design and residual design of the unique 4-$(23,7,1)$ design $(V^+,\mathcal B)$.\\
{\bf{case 1: $r_1=6,r_2=7$.}}\\
Assume $Y=Y_6 \cup Y_7$ is a tight relative 4-design and $(V,\mathcal B_6)$ corresponds to $Y_6$.
It is to prove that $(V,\mathcal B_7)$ corresponds to $Y_{r_2}$.
Suppose $(V,\mathcal B_7^{\prime})$ is a different 3-$(22,7,4)$ design cooresponding to $Y_7^{\prime}$ such that $Y=Y_6 \cup Y_7^{\prime}$.
We should remark that $(V,\mathcal B_7^{\prime})$ may be isomorphic to $(V,\mathcal B_7)$.
Using Lemma \ref{lem:extension}, we have two different 4-$(23,7,1)$ designs $(V^+,\mathcal B):=(V^+,\mathcal B_6^+ \cup \mathcal B_7)$ and $(V^+,\mathcal B^{\prime}):=(V^+,\mathcal B_6^+ \cup \mathcal B_7^{\prime})$.\\
(i). If $(V,\mathcal B_7^{\prime})$ and $(V,\mathcal B_7)$ are non-isomorphic, then $(V^+,\mathcal B)$ and $(V^+,\mathcal B^{\prime})$ are non-isomorphic. This contradicts the uniqueness of 4-$(23,7,1)$ design.\\
(ii). If $(V,\mathcal B_7^{\prime})$ is isomorphic to $(V,\mathcal B_7)$ (but different), then $\exists$ a permutation $\sigma \in S_{22}$ such that $\sigma(\mathcal B_7)=\mathcal B_7^{\prime}$, where $S_{22}$ is the symmetric group on the vertex set $V$.
This implies that $(V^+,\mathcal B)$ and $(V^+,\mathcal B^{\prime})$ are not isomorphic.
Otherwise there are at least two distinct points $p_1, p_2 \in V$ which are not fixed by $\sigma$, such that $p_1, p_2$ appear in the same ($=21$) blocks in $\mathcal B_6$.
It is impossible because $\lambda_2=5$ for any 3-$(22,6,1)$ design.
This again gives a contradiction!\\
{\bf{case 2: $r_1=6,r_2=15$.}}\\
Assume $Y=Y_6 \cup Y_{15}$ is a tight relative 4-design and and $(V,\mathcal B_6)$ corresponds to $Y_6$.
We have verified that $Y_6 \cup Y_7^c$ is a tight relative 4-design.
It is enough to prove that $(V,\mathcal B_7^c)$ corresponds to $Y_{15}$.
Suppose $(V,\mathcal B_{15}^{\prime})$ is a different 3-$(22,15,52)$ design corresponding to $Y_{15}$.
With the similar argument in case 1, we obtain two non-isomorphic 4-$(23,7,1)$ designs $(V^+,\mathcal B_6^+ \cup \mathcal B_7^c)$ and $(V^+,\mathcal B_6^+ \cup \mathcal B_{15}^{\prime})$.
This again gives a contradiction.\\
The other two cases can also be proved using the uniqueness of 4-$(23,15,52)$ design.\\

\noindent
Finally, we exclude the tight relative 4-designs with non-constant weight.
It is not difficult to check that the types of $\lambda_4$-sequence are $(6160*0,1155*1)$, $(6160*1,1155*0)$, $(6160*33,1155*32)$, $(6160*19,1155*20)$ corresponding to $r=6,7,15,16$.
Then the list below implies the non-existence of tight relative 4-designs with non-constant weight for $n=22$.

\begin{spacing}{1.1}
$$
{\footnotesize
\begin{array}{ccc|cc|c|c}
n & r_1 & r_2 &\lambda_3^{r_1} & \lambda_3^{r_2}  & \frac{w_{r_2}}{w_{r_1}} & (\lambda_4^{(r_1)},\lambda_4^{(r_2)}) \\ \hline
 22 & 6 & 15 & 1 & 52 & \frac{1}{20} & (0,36),(1,16)\\ \hline
 22 & 7 & 16 & 4 & 28 &\frac{4}{23} & (0,24),(4,1)\\ \hline
   &   &    &    &    &\frac{2}{21} & (0,28),(2,7)\\ \hline
  22 & 15 & 16 & 52 &28 & 39 & (0,20),(39,19) \\ \hline
& & & & & \frac{24}{5} & (0,26), (24,21),(48,16)\\ \hline
& & & & & \frac{26}{7} & (0,28), (26,21), (52,14)\\ \hline
& & & & & \frac{21}{2} & (3,22), (24,20), (45,18)\\ \hline
& & & & & \frac{27}{8} & (3,28), (30,20)\\ \hline
& & & & & \frac{23}{4} & (5,24), (28,20), (51,16)\\ \hline
& & & & & \frac{31}{12} & (10,28), (41,16)\\ \hline
& & & & & \frac{22}{3} & (12,22), (34,19)\\ \hline
& & & & & \frac{33}{14} & (12,28), (45,14)\\ \hline
& & & & & \frac{29}{10} & (13,26), (42,16)\\ \hline
& & & & & 20 & (16,20), (36,19)\\ \hline
& & & & & \frac{32}{13} & (16,26), (48,13)\\ \hline
& & & & & \frac{25}{6} & (21,22), (46,16)\\ \hline
& & & & & \frac{28}{9} & (24,22), (52,13)\\ \hline
& & & & & \frac{5}{24} & (31,28), (36,4)\\ \hline
& & & & & \frac{4}{23} & (32,24), (36,1)\\ \hline
& & & & & \frac{2}{21} & (32,28), (34,7)\\ \hline
\end{array}}
$$
\end{spacing}
This completes the proof. \hfill\qed

\section{Concluding Remarks}

Here we collect some open problems we want to study in the research direction described in the present paper.

Firstly, let us recall the three problems (1), (2), (3) mentioned at the end of Section \ref{sec:main}.
We would like to add the following problems.

(i) So far, we do not know any example of tight relative 3-designs on two shells in $H(n,2)$ with non-constant weight function, more precisely constant on each shell.
The first several open cases are $n=56,$ then $n=66$ (Case 2) and $n=61$ (Case 4).
We bet that there may exist such tight relative 3-designs on two shells in $H(n,2)$ with non-constant weight, but so far we have difficulty in finding one.

(ii) So far, we do not know any counterexample to Conjecture 1 and Conjecture 2. The first open cases seem to be $n=36$ and then $n=40$.
Although we formulated Conjecture 1 and Conjecture 2, we think these are important as working hypothesis, and we would not be surprised even if a counterexample could be found.
It would be theoretically very interesting whether these conjectures hold or not.

(iii) Conjecture 2 obviously does not hold, if the condition that the two 2-$(v,k,\lambda)$ designs are symmetric is dropped.
For example, if you take any two non-isomorphic 3-designs, then regarding these designs as 2-designs, they are not isomorphic as 2-designs, but obviously have the same $\lambda_3$-sequence.
So, it would be interesting to try to find two non-isomorphic 2-designs close to symmetric designs with different $\lambda_3$-sequences.
It would be interesting for which family of symmetric designs, or non-symmetric designs, whether the property mentioned in Conjecture 2 holds or not.

(iv) It would be interesting how much the methods used in the present paper could be generalized for the study of tight relative $t$-designs on other Q-polynomial association schemes.
In particular, it would be interesting to know how much theorems similar to our Theorem \ref{theo:3.3} (as well as Kageyama's theorem) hold for other association schemes.
First test cases would be non-binary Hamming association schemes $H(n,q)$ and Johnson association schemes $J(v,k).$

(v) We would like to repeat our belief that the classification problem of tight relative $t$-designs is interesting problem as the classification problem of tight $t$-designs is interesting.

\section*{Acknowledgments}
 We are very grateful to many people for their help in completing this work.
 We thank Andries Brouwer and Akihiro Munemasa for their help for us to understand Driessen's work on certain 3-designs.
 In particular, Andries Brouwer kindly sent us the Ph. D thesis of Driessen \cite{Driessen-1978}.
 We are very grateful to Ted Spence for sending us the incidence matrices of all the $1106$ number of $2$-$(23,11,5)$ designs and all the $151$ number of $2$-$(31,10,3)$ designs, which were very useful in our considerations of these cases.
 We also very much benefited from the information presented in his home page.
 We are very thankful to Akihide Hanaki and Hadi Kharaghani answering our questions on $2$-$(31,15,7)$ designs.
 We are extremely grateful for Brendan McKay for showing us how to obtain all $2$-$(31,15,7)$ designs from the list of all Hadamard matrices of order $32$ by Kharaghani and Tayfeh-Rezaie.
 Moreover, Brendan McKay kindly determined all the $\lambda_3$-sequences with length at most $4$, among all the $2$-$(31,15,7)$ designs for us.
 They were extremely useful for our research presented in this paper.
 We also thank Zongchen Chen and Da Zhao, undergraduate students of Shanghai Jiao Tong University.
 They made a decisive contribution to the problem presented in this paper.
 Their results will be published independently, but we included the explanation of some of their results in this paper.
 The first author (Eiichi Bannai) is supported in part by NSFC grant No. 11271257.

\noindent
Eiichi Bannai: Department of Mathematics, Shanghai Jiao Tong University,
800 Dongchuan Road, Shanghai, 200240, China\\
e-mail: bannai@sjtu.edu.cn
\\

 \noindent
Etsuko Bannai: Misakigaoka 2-8-21, Itoshima-shi, Fukuoka, 819-1136, Japan\\
e-mail: et-ban@rc4.so-net.ne.jp
\\

\noindent
Yan Zhu: Department of Mathematics, Shanghai Jiao Tong University,
800 Dongchuan Road, Shanghai, 200240, China\\
e-mail: zhuyan@sjtu.edu.cn
  \end{document}